\documentclass[preprint,9pt]{elsarticle}
\usepackage{amssymb}
\usepackage{amsmath}
\usepackage{amsthm}
\usepackage{hyperref}
\usepackage{amsfonts}
\usepackage{color}
\usepackage{epsfig}
\usepackage{graphics,subfigure}                 
\usepackage{graphicx}
\usepackage{float}
\usepackage{epsf}
\usepackage{mathrsfs}
\allowdisplaybreaks[0]
\biboptions{numbers,sort&compress}

\newtheoremstyle{thmm}{1.5ex plus 1ex minus .2ex}{1.5ex plus 1ex minus
.2ex}{\rmfamily}{}{\bfseries}{}{1em}{} \theoremstyle{thmm}
\newtheorem{theorem}{Theorem}[section]
\newtheorem{lemma}{Lemma}[section]

\newtheorem{remark}{Remark}[section]
\newtheorem{assumption}{Assuption}[section]
\newtheorem{example}{Example}[section]

\allowdisplaybreaks
\def\refe#1{(\ref{#1})}

\def\u{\mathbf{u}}
\def\b{\mathbf{b}}
\def\a{\mathbf{a}}
\def\r{\mathbf{r}}
\def\d{\mathbf{d}}
\def\v{\mathbf{v}}
\def\c{\mathbf{c}}
\def\w{\mathbf{w}}
\def\q{\mathbf{q}}
\def\f{\mathbf{f}}
\def\m{\mathbf{m}}
\def\g{\mathbf{g}}
\def\x{\mathbf{x}}
\def\n{\mathbf{n}}
\def\0{\mathbf{0}}
\def\H{\mathbf{H}}
\def\V{\mathbf{V}}
\def\M{\mathbf{M}}
\def\R{\mathbb{R}}
\def\L{\mathbf{L}}
\def\Z{\mathbf{Z}}
\def\B{\mathbf{B}}
\def\W{\mathbf{W}}
\def\I{\mathrm{I}}
\def\P{\mathbb{P}}
\def\G{\mathcal{G}}
\def\Th{\mathcal{T}_h}
\def\Vh{\V_h}
\def\Mh{\M_h}
\def\Zh{\Z_h}
\def\uh{\u_h}
\def\bh{\b_h}
\def\dh{\d_h}
\def\ah{\a_h}
\def\rrh{\r_h}
\def\wh{\w_h}
\def\vh{\v_h}
\def\ch{\c_h}
\def\fh{\f_h}
\def\gh{\g_h}
\def\hu{\overline{\u}}
\def\hb{\overline{\b}}
\def\hq{\overline{q}}
\def\he{\overline{e}}

\def\huh{\hu_h}
\def\hbh{\hb_h}
\def\hqh{\hq_h}
\def\qh{q_h}
\def\rh{r_h}
\def\curl{\mathrm{curl\,}}
\def\div{\mathrm{div\,}}
\def\pt{\partial}
\def\O{\Omega}
\def\din{\mathrm{in\;}}
\def\ds{\displaystyle}
\def\dta{\Delta}

\def\error{\mathrm{error}}

\begin{document}
\begin{frontmatter}
\title{\bf  A virtual element method with IMEX-SAV scheme for the
incompressible magnetohydrodynamics equations
\thanks{.}}

\author[]{Xiaojing Dong}
\ead{dongxiaojing99@xtu.edu.cn}

\author[]{Yunqing Huang}
\ead{huangyq@xtu.edu.cn}

\author[]{Tianwen Wang \corref{mycorrespondingauthor}}
\cortext[mycorrespondingauthor]{Corresponding author}
\ead{wangtianwen@smail.xtu.edu.cn}

\address{Hunan Key Laboratory for Computation and Simulation in Science and Engineering, Key Laboratory of Intelligent Computing $\&$ Information Processing of Ministry of Education, School of Mathematics and Computational Science, Xiangtan University, Xiangtan, Hunan, 411105, P.R. China}

\begin{abstract}
This paper proposes a virtual element method (VEM) combined with a second-order implicit-explicit scheme based on the scalar auxiliary variable (SAV) method for the incompressible magnetohydrodynamics (MHD) equations. We employ the BDF2 scheme for time discretization and a conservative VEM for spatial discretization, in which the mass conservation in the velocity field is kept by taking advantage of the virtual element method’s adaptability and its divergence-free characteristics. In our scheme, the nonlinear terms are handled explicitly using the SAV method, and the magnetic field is decoupled from the velocity and pressure. This decoupling only requires solving a sequence of linear systems with constant coefficient at each time step. The stability estimate of the fully discrete scheme is developed, demonstrating the scheme is unconditionally stable. Moreover, rigorous error estimates for the  velocity and magnetic field are provided. Finally, numerical experiments are presented to verify the valid of theoretical analysis.

\end{abstract}
\end{frontmatter}
\section{Introduction}
We consider the following two-dimensional incompressible magnetohydrodynamics (MHD) equations:
\begin{align} \label{MHDeq}
\left\{
\begin{aligned}
\u_t - \nu \Delta  \u + (\nabla \u) \u + \nabla p - \mu \curl \b \times \b &= \f \quad \din \O, \\
\mu \b_t + \sigma^{-1}\curl\curl\b - \mu \curl(\u\times\b) &= \g \quad \din\O, \\
\div \u = 0, \quad \div \b &= 0  \quad  \din \O, 
\end{aligned}
\right.
\end{align}
for all times $t\in (0,T]$, where $\u$ is the velocity, $\b$ the magnetic field, $p$ the pressure, $\nu$ the kinematic viscosity, $\sigma$ the magnetic permeability and $\mu$ the electric conductivity. Functions $\f$ and $\g$ represent the external body forces. The equations \refe{MHDeq} satisfy the following initial boundary conditions
\begin{align*}
& \u(\x,0) = \u_0(\x),\quad \b(\x,0) = \b_0(\x), \quad \forall\, \x \in \O,\\
& \u = \0, \quad \b\cdot \n = 0,\quad \n \times \curl \b = \0,  \quad \forall\, \x \in \pt \O,\, t\in [0,T],
\end{align*}
where $\n$ is the unit outer normal of $\pt \O$ and $\div \u_0(\x) = 0, \div \b_0(\x)=0$.

The incompressible MHD equations, coupling Naiver-Stokes equations and Maxwell's equations, model the interaction of a viscous, imcompressible, electrically conducting fluid and an applied magnetic field. The incompressible MHD equations play an important role in astronomy and geophysics as well as engineering problems, such as the propagation of radio waves in the ionosphere, cooling of liquid metal in nuclear reactors, confinement for controlled thermonuclear fusion and magnetic fluid motors, see \cite{MHD-introduce,MHD-book}. Due to the nonlinearity of convection term and coupling terms, the numerical approximation of MHD equations is a challenging work. In recent decades, there has been a lot of research on the numerical solution of MHD equations, see \cite{MHD-existence,He-MHD-Euler,Dong-MHD-iterative,Dong-CN-MHD,Dong-SAV-MHD,Tang-MHD,Xu-projection-MHD,Li-SAV-MHD,Li-MHD-divergence-free,VEM-3D-MHD,VEM-2D-MHD} and the references therein. As everyone knows, for unsteady MHD equations, a fully implicit scheme offers unconditional stability, but needs to solve a nonlinear system at each time step. In contrast, while an explicit scheme is simpler to implement, it generally only offers conditional stability. Consequently, many researchers are interested in semi-implicit schemes for MHD equations, as they strike a balance between stability and computational efficiency. In \cite{He-MHD-Euler}, a first-order Euler semi-implicit scheme based on the mixed finite element method (FEM) was developed. A Crank–Nicolson extrapolation scheme was studied in \cite{Dong-CN-MHD}, almost unconditionally stable and the optimal error estimates were provided. A fully divergence-free FEM with Crank–Nicolson scheme was introduced in \cite{Li-MHD-divergence-free}. The literature \cite{Xu-projection-MHD} proposed and analyzed a Crank–Nicolson finite element projection method.  

The scalar auxiliary variable (SAV) method for the gradient flow was proposed in  \cite{Shen-SAV}. The SAV method is attractive for the numerical approximation of nonlinear systems. Since the method allows to design unconditionally stable schemes, but only requires solving decoupled linear systems with constant coefficients at each time step. Reference \cite{Li-SAV-NS} presented the SAV-pressure correction methods for Navier-Stokes. In \cite{Hou-SAV-NS}, a $\H^1$-conforming divergence-free FEM with ﬁrst-order implicit-explicit (IMEX) SAV scheme for Navier-Stokes equations at high Reynolds number was analyzed. The development of first- and second-order IMEX SAV schemes for MHD equations, along with error estimates for the first-order SAV semi-discrete scheme, was detailed in \cite{Li-SAV-MHD}. Additionally, reference \cite{Dong-SAV-MHD} investigated two rotational pressure-correction SAV schemes specifically for MHD equations. 

Virtual element method (VEM),  first introduced in \cite{VEM-2013},  is a novel numerical method, which allows for numerical simulation on greatly general polygonal meshes (including convex or nonconvex elements). The computer implementation of VEM for the linear elliptic problem  was presented in \cite{VEM-2014}. The Stoke-like virtual element space was constructed and applied to Stokes equations \cite{DivFree-Stokes-CVEM}, with enhanced version later used for Navier-Stokes equations in \cite{DivFree-NS-CVEM}. Due to the divergence-free property of the virtual element space within the associated discrete kernel, this space is widely used for incompressible problems. In \cite{stream-VEM-NS}, a VEM for the time dependent Navier–Stokes equations in stream-function formulation was proposed. The fully divergence-free VEMs with implicit Euler scheme for resistive MHD equations were developed in \cite{VEM-3D-MHD,VEM-2D-MHD}. 

The aim of this paper is to present a conservative VEM with a BDF2 IMEX scheme based on the SAV method for the incompressible MHD equations. For the time and space discretizations, the BDF2 scheme and VEM are adopted, respectively. We employ the enhanced Stokes-like virtual element to approximate the velocity, discontinuous piecewise polynomials for the pressure, $\H^1$-conforming virtual element \cite{VEM-equivalent} for the magnetic field. The nonlinear terms are treated explicitly based SAV method, to ensure the high efficiency and unconditional stability of scheme. In our method, the magnetic field is decoupled from the velocity and pressure,  thus, we only need to solve the Stokes and Poisson subproblems with constant coefficients at each time step. Moreover, the stability estimate and the optimal error estimates of the fully discrete scheme are provided.


The rest of this paper is organized as follows. In Section 2, we describe some basic notations and the weak formulation of the problems \refe{MHDeq}. In Section 3, the virtual element approximation spaces of the velocity,  pressure and magnetic field are introduced. We give the definitions of discrete bilinear and trilinear forms and load terms, and discuss the related properties. In Section 4, we present the fully discrete BDF2 IMEX SAV scheme and establish the stability estimate. Section 5 is devoted to providing the rigorous error estimates. Finally, Section 6 shows some numerical experiments to validate the theoretical results.

\setcounter{equation}{0}
\section{Preliminaries}
Consider a polygonal domain $D \in \R^2$, $W^{s,l}(D)$ is the scalar Sobolev space equipped with norm $\|\cdot\|_{s,l,D}$ and semi-norm $|\cdot|_{s,l,D}$ for $s \geq 0$ and $l\geq 1$. If $s>0,l=2$, we denote $W^{s,l}(D)$ by $H^s(D)$ equipped with norm $\|\cdot\|_{s,D}$ and semi-norm $|\cdot|_{s,D}$. If $s=0$, the $W^{s,l}(D)$ is the usual Lebesgue space $L^l(D)$ with norm $\|\cdot\|_{L^l(D)}$. The inner product and the norm in $L^2(D)$ are denoted by $(\cdot,\cdot)_{D}$ and $\|\cdot\|_{0,D}$, respectively. When $D$ is the whole computational domain $\O$, the subscript $\O$ will be omitted in inner product, norm and semi-norm. $\mathbf{K}$ stands for the corresponding vectorial version of a generic scalar $K$, examples of this are $\W^{s,l}(\O) = [W^{s,l}(\O)]^2$ and $\L^2(\O)= [L^2(\O)]^2$. Let $t$ represent the time variable taking values in $[0,T]$, $X$ is a real Banach space with the norm $\|\cdot\|_X$. For the function space $L^l(0,T;X)$, $ 1 \leq l \leq \infty$, the norm is denoted as
\begin{align*}
\|\cdot\|_{L^l(0,T;X)} = 
\begin{cases}
\left(\ds\int_0^T \|\cdot\|_X^l \,dt \right)^{\frac{1}{l}},\qquad 1 \leq l < \infty, \\
\underset{t\in [0,T]}{\mathrm{ess}\, \sup} \|\cdot\|_X, \qquad l = \infty.
\end{cases}
\end{align*}  
Throughout the paper, $C$ is a generic positive constant independent of mesh size.

Before showing weak formulation, the following spaces are introduced
\begin{align*}
\V &= \H_0^1(\O) = \{ \v \in \H^1(\O):\; \v|_{\pt \O} = \0\}, \qquad
\M = \H_n^1(\O) = \{ \c \in \H^1(\O):\; \c\cdot \n|_{\pt \O} = 0\},\\
\Z &= \{\v \in \H_0^1(\O):\; \div \v = 0 \}, \qquad
\B = \{\c \in \H_n^1(\O): \div \c = 0 \} , \\
Q &=  L_0^2(\O) = \{  q \in  L^2(\O):\;  (q,1) = 0\} .
\end{align*}
The weak formulation of \refe{MHDeq} can be rewritten as follows: for almost every time $t\in (0,T]$, find $(\u(t),\b(t),p(t))\in \V\times\M\times Q$ such that
\begin{align}
\begin{aligned}
m_0(\u_t,\v) + a_0(\u,\v) - c_0(\u,\u,\v) - c_1(\b,\b,\v)+  d(\v,p)&= (\f,\v), \quad \forall\, \v \in \V,  \\
m_1(\b_t,\c) + a_1(\b,\c) + c_1(\c,\b,\u) &= (\g,\c), \quad \forall\, \c \in \M,\\
d(\u,r) &= 0, \quad \forall\, r \in Q,
\end{aligned} \label{weakMHDeq}
\end{align} 
where, $\forall\, \u,\v \in \V,\; \w \in \Z,\; \b,\c \in \M,\; r \in Q$,
\begin{align*}
&m_0(\u,\v) = (\u,\v), \quad m_1(\b,\c) = \mu (\b,\c), \quad a_0(\u,\v) = \nu (\nabla \u,\nabla \v), \\
&a_1(\b,\c) = \sigma^{-1}(\curl \b, \curl \c) + \sigma^{-1}(\div \b, \div \c), \quad c_0(\w,\u,\v) = ((\nabla \v)\w,\u), \\
&c_1(\c,\b,\u)= \mu(\curl \c\times \b,\u) ,\quad d(\v,r) = -(\div \v,r). 
\end{align*}
There holds
\begin{align}
\curl (\b\times\d) = (\d \cdot \nabla)\b - (\b \cdot \nabla)\d,\quad \forall\, \b,\d \in \B. \label{curl-eq1}
\end{align}

\setcounter{equation}{0}
\section{Virtual element method}
\subsection{Basis setting}
Let $\Th$ be a decomposition of $\O$ into non-overlapping polygons $E$, with mesh size $h = \max_{E \in \Th} h_E$, where $h_E$ denotes the diameter of the element $E$. We denote each edge of $\pt E$ by $e$ and the measure of the element $E$ by $|E|$. Moreover, each element $E$ of $\Th$ satisfies the following assumptions: 
\begin{flushleft}
  (\textbf{A1}) $E$ is star-shaped with respect to a ball $B_E$ of radius $\ge \varrho h_E$,  \\
  (\textbf{A2}) the distance between any two vertexes of $ E$ is $\geq \varrho h_E$,
\end{flushleft} 
where $\varrho$ is a positive constant. For each integer $k \geq 0$, we define the polynomial spaces
\begin{itemize}
\item $\P_k(E)$ stands for the set of polynomials on element $E$ of degrees $\leq k$, specially, $\P_{-2}(E) = \P_{-1}(E) = \emptyset $,
\item $\G_k(E) = \nabla \P_{k+1}(E),\; \G_k^{\oplus}(E) = \x^{\perp}\P_{k-1}(E)$, where $\x^{\perp}=  (-x_2,x_1)^{\mathrm{T}}$,
\end{itemize} 
where $\mathrm{T}$ represents the standard transpose. In order to define computable discrete linear forms, we introduce some useful projections as follows:
\begin{itemize}
   \item $H^1$-semi-norm projection $\Pi_k^{\nabla,E}$ : $\H^1(E) \rightarrow [\mathbb{P}_k(E)]^2$, defined by
   \begin{align}
   \begin{cases}
    \int_E \nabla \Pi_k^{\nabla,E} \v : \nabla \q_k \,\mathrm{d} E = \int_E \nabla \v : \nabla \q_k \,\mathrm{d}E, \quad \forall\, \q_k \in \left[ \P_k(E) \right]^2, \;\v \in \H^1(E),\\
    \int_{\pt E} \Pi_k^{\nabla,E} \v \,\mathrm{d} \pt E = \int_{\pt E} \v \;\mathrm{d}\pt E ,\quad \forall \,\v \in \H^1(E).
   \end{cases}   
   \end{align}
   \item $L^2$-projection $\Pi_k^{0,E}$ : $\L^2(E) \rightarrow  \left[\P_k(E)\right]^2 $, defined by
   \begin{align}
   \int_E \Pi_k^{0,E} \v \cdot \q_k \,\mathrm{d}E = \int_E\v \cdot \q_k \,\mathrm{d}E ,  \quad \forall\, \q_k \in \left[ \P_k(E) \right]^2,\, \v \in \L^2(E). 
   \end{align}
\end{itemize}
The above $L^2$-projection of  vector functions can be extended to the scale and tensor functions.

\subsection{Virtual element spaces}
\subsubsection{The discrete velocity and pressure spaces}
We use the enhanced Stokes-like virtual element \cite{DivFree-NS-CVEM,VEM-2D-MHD} of order $k \geq 1$, which is divergence-free in associate discrete kernel, to approximate the velocity. In order to introduce the element, we first define the following local virtual element
\begin{align*}
\widetilde{\V}_h(E) = \Big\{\v \in \H^1(E)  :  \v|_{\partial E} \in C^0(\partial E) ,\; \v|_e \in [\mathbb{P}_{\hat{k}}(e)]^2, \;\; \forall e \in \partial E, \;\div \v \in \mathbb{P}_{k-1}(E),\\
-\Delta \v - \nabla \phi \in \G_k^{\oplus}(E) \quad \mathrm{for} \;\mathrm{some} \;\phi \in L^2(E)\Big\},
\end{align*} 
where $\hat{k} = \max\{2,k\}$. The local enhanced Stokes-like virtual element space is defined as
\begin{align*}
\V_h(E) = \left\{\v \in \widetilde{\V}_h(E): (\v - \Pi_k^{\nabla ,E}\v,\q)_E=0, \quad \forall\, \q \in \G_k^{\oplus}(E)/\G_{k-2}^{\oplus}(E) \right\}, 
\end{align*}
with the degrees of freedom
\begin{align*}
& \bullet \quad \mathrm{the} \; \mathrm{values} \;  \mathrm{of} \; \v \; \mathrm{at} \; \mathrm{vertices}\;  \mathrm{of}\; E,      \\
& \bullet \quad  \mbox{the}\; \mbox{values}\; \mbox{of}\; \v \;\mbox{at}\; \hat{k}-1\; \mbox{distinct}\; \mbox{points}\; \mbox{of}\; \mbox{every}\; \mbox{edge}\; e \;\mbox{of}\; \pt E, \\
& \bullet \quad \frac{1}{|E|}\int_E  \v \cdot  \q_{k-2} \,\mathrm{d}E ,\quad \forall\,  \q_{k-2} \in \G_{k-2}^{\oplus}(E),  \\
& \bullet \quad \frac{h_E}{|E|}\int_E (\div \v)\, q_{k-1} \,\mathrm{d}E , \quad \forall\, q_{k-1}\in \P_{k-1}(E)/\mathbb{R}. 
\end{align*}
As everyone knows, the above degrees of freedom allow us to compute exactly the projections $\Pi_k^{\nabla,E}\circ \V_h(E)$, $\Pi_k^{0,E}\circ \V_h(E)$ and $\Pi_{k-1}^{0,E} \circ \nabla \V_h(E)$.

The global space is defined by
\begin{align}
\V_h = \left\{\v \in \H_0^1(\O) : \v |_E \in \V_h(E), \quad \forall\, E \in \Th\right\}. \label{globalDFCVEM}
\end{align}

For the discrete pressure space, we choose the discontinuous piecewise polynomials
\begin{align}
Q_h = \left\{ r \in L_0^2(\O): \,r|_E \in \P_{k-1}(E) ,\quad \forall \,E \in \Th \right\}.
\end{align} 
Obviously, we can observe that $\div \V_h \in Q_h$, which is crucial to divergence-free of the discrete velocity. We define the discrete kernel space
\begin{align}
\Z_h = \{ \vh \in \V_h: d(\vh,r_h) = 0, \quad \forall\, r_h \in Q_h \}.
\end{align}
From \cite{DivFree-NS-CVEM,VEM-2D-MHD}, we know that the spaces $\V_h$ and $Q_h$ form an inf-sup stable virtual element  pair on polygonal meshes.
\subsubsection{The discrete magnetic field space}
For computability of $L^2$-projection, the enhanced node virtual element \cite{VEM-equivalent}, which is modified based on the original node space \cite{VEM-2013}, is used to approximate the magnetic field. First, we introduce
\begin{align*}
\tilde{M}_h(E) = \left\{ m \in H^1(E) :\, \Delta m \in \mathbb{P}_{k}(E),\,   m|_{\partial E} \in C^0(\partial E) ,\,m|_e \in \mathbb{P}_{k}(e) ,\quad \forall \, e \in \partial E   \right\} .
\end{align*}
Then, the local enhanced space of order $k \geq 1$ is defined by 
\begin{align*}
\M_h(E) = \left\{ \m\in [\tilde{M}_h(E)]^2:\; (\m - \Pi_k^{\nabla,E} \m,\q)_E=0, \quad \forall \,\q \in [\mathbb{P}_k(E)]^2/[\mathbb{P}_{k-2}(E)]^2 \right\}, 
\end{align*}
with the degrees of freedom 
\begin{align*}
& \bullet \quad \mathrm{the} \; \mathrm{values} \;  \mathrm{of} \; \m \; \mathrm{at} \; \mathrm{vertices}\;  \mathrm{of}\; E,      \\
& \bullet \quad  \mbox{the}\; \mbox{values}\; \mbox{of}\; \m \;\mbox{at}\; k-1\; \mbox{distinct}\; \mbox{points}\; \mbox{of}\; \mbox{every}\; \mbox{edge}\; e \;\mbox{of}\;  \pt E, \\
& \bullet \quad \frac{1}{|E|} \int_E  \m \cdot  \q_{k-2} \,\mathrm{d}E ,\quad \forall \, \q_{k-2} \in [\P_{k-2}(E)]^2.  
\end{align*}
According to integration by parts and the above degrees of freedom, we see that these projections $\Pi_k^{\nabla,E} \circ \M_h(E)$, $\Pi_k^{0,E} \circ \M_h(E), \Pi_{k-1}^{0,E}\circ \curl \M_h(E)$ and $\Pi_{k-1}^{0,E}\circ \div \M_h(E)$ are computable. From the properties of $L^2$-projection, we observe that for all $\ch \in \M_h$, there hold
\begin{align}
\begin{aligned}
\|(\I - \Pi_{k-1}^{0,E})\curl \ch\|_{0,E} &\leq \|\curl (\I - \Pi_k^{\nabla,E})\ch\|_{0,E} \leq \sqrt{2}\|\nabla (\I - \Pi_k^{\nabla,E})\ch\|_{0,E},\\
\|(\I - \Pi_{k-1}^{0,E})\div \ch\|_{0,E} &\leq \|\div (\I - \Pi_k^{\nabla,E})\ch\|_{0,E} \leq \sqrt{2}\|\nabla (\I - \Pi_k^{\nabla,E})\ch\|_{0,E}. 
\end{aligned} \label{projection-eq}
\end{align}

The global space is defined as
\begin{align}
\M_h = \left\{\m \in \H_n^1(\O):\; \m|_E \in \M_h(E)  ,\quad \forall\, E \in \Th   \right\}.\label{globalCVEM}
\end{align}

\subsection{Discrete bilinear and trilinear forms and load terms approximation}
We use $H^1$-semi-norm projection and $L^2$-projection to define the computable discrete local bilinear forms, for all $\uh,\vh \in \V_h(E),\bh,\ch \in \M_h(E) $,
\begin{align}
m_{0h}^E(\uh,\vh) &= m_{0}^E(\Pi_k^{0,E}\uh,\Pi_k^{0,E}\vh) + |E|S_0^E((\I - \Pi_k^{0,E})\uh,(\I - \Pi_k^{0,E})\vh),\label{bilinear-eq1}\\
m_{1h}^E(\bh,\ch) &= m_{1}^E(\Pi_k^{0,E}\bh,\Pi_k^{0,E}\ch) + \mu|E|S_1^E((\I - \Pi_k^{0,E})\bh,(\I - \Pi_k^{0,E})\ch),\label{bilinear-eq2}\\
a_{0h}^E(\uh,\vh) &= a_{0}^E(\Pi_k^{\nabla,E}\uh,\Pi_k^{\nabla,E}\vh)+ \nu S_0^E((\I - \Pi_k^{\nabla,E})\uh,(\I - \Pi_k^{\nabla,E})\vh),\label{bilinear-eq3}\\ 
a_{1h}^E(\bh,\ch) &= \sigma^{-1} (\Pi_{k-1}^{0,E}\curl\bh,\Pi_{k-1}^{0,E}\curl\ch)_E + \sigma^{-1}(\Pi_{k-1}^{0,E}\div\bh,\Pi_{k-1}^{0,E}\div\ch)_E \notag\\
 &\quad + \sigma^{-1} S_1^E((\I - \Pi_k^{\nabla,E})\bh,(\I - \Pi_k^{\nabla,E})\ch), \label{bilinear-eq4}
\end{align}
where $S_0^E(\cdot,\cdot)$ and $S_1^E(\cdot,\cdot)$ are symmetric positive definite bilinear forms such that
\begin{align*}
\eta_{0\ast}m_0^E(\vh,\vh) \leq |E|S_0^E(\vh,\vh) \leq \eta_0^{\ast}m_0^E(\vh,\vh), \quad \forall  \vh \in \V_h(E)\cap \mbox{Ker}(\Pi_k^{0,E}),  \\
\eta_{1\ast}m_1^E(\ch,\ch) \leq \mu|E|S_1^E(\ch,\ch) \leq \eta_1^{\ast}m_1^E(\ch,\ch), \quad \forall  \ch \in \M_h(E)\cap \mbox{Ker}(\Pi_k^{0,E}),\\
\eta_{2\ast}a_0^E(\vh,\vh) \leq \nu S_0^E(\vh,\vh) \leq \eta_2^{\ast}a_0^E(\vh,\vh), \quad \forall  \vh \in \V_h(E)\cap \mbox{Ker}(\Pi_k^{\nabla,E}),  \\
\eta_{3\ast}\|\nabla \ch\|_{0,E} \leq S_1^E(\ch,\ch) \leq \eta_3^{\ast}\|\nabla \ch\|_{0,E}, \quad \forall  \ch \in \M_h(E)\cap \mbox{Ker}(\Pi_k^{\nabla,E}).
\end{align*}
for positive constants $\eta_i^{*}, \eta_{i*} (i=0,1,2,3)$ independent of $h$ and $E$. For the stabilizing terms $S_0^E(\cdot,\cdot)$ and $S_1^E(\cdot,\cdot)$, we follow the classical choice
\begin{align*}
S_0^E(\uh,\vh) = \sum_{i=1}^{\mathrm{dim}(\V_h(E))}\mathrm{dof}_i^{\V_h(E)}(\uh)\mathrm{dof}_i^{\V_h(E)}(\vh),\\
S_1^E(\bh,\ch) = \sum_{i=1}^{\mathrm{dim}(\M_h(E))}\mathrm{dof}_i^{\M_h(E)}(\bh)\mathrm{dof}_i^{\M_h(E)}(\ch),
\end{align*}
where the operator $\mathrm{dof}_i^{\W}(\v)$ is the $i$ local degree of freedom of smooth enough function $\v$ to the local space $\W$, here, $\W = \V_h(E)$ or $\M_h(E)$. Other choices for the stabilizing terms can be found in \cite{stability-MMM,interpolation-Stokes,stability-ill,chen-huang-2018}. From \cite{VEM-equivalent,VEM-2013,DivFree-NS-CVEM,chen-huang-2018,CVEM-elliptic,CVEM-NCVEM-elliptic,interpolation-Stokes} and \refe{projection-eq}, we know that the above bilinear forms satisfy the following properties.
\begin{lemma} \label{stalility-consistency}
The discrete local bilinear forms  \refe{bilinear-eq1}-\refe{bilinear-eq4} satisfy the following properties 
\begin{itemize}
\item $k$-consistency: for all $\vh \in \V_h(E), \ch \in \M_h(E),\q_k \in \left[\P_k(E) \right]^2$ and $E \in \Th$, we have
\begin{align*}
m_{0h}^E(\vh,\q_k) &= m_0^E(\vh,\q_k),\quad m_{1h}^E(\ch,\q_k) = m_1^E(\ch,\q_k),\\
a_{0h}^E(\vh,\q_k) &= a_0^E(\vh,\q_k),\quad a_{1h}^E(\ch,\q_k) = a_1^E(\ch,\q_k).
\end{align*}
\item stability: there exist positive constants $\kappa_i^{*}, \kappa_{i*}, i=0,1,2,3$, independent of $h$ and $E$ satisfying
\begin{align*}
&\kappa_{0*} m_{0}^E(\vh,\vh)\leq m_{0h}^E(\vh,\vh) \leq \kappa_0^* m_{0}^E(\vh,\vh), \quad \forall\,\vh \in \V_h(E),\\
&\kappa_{1*} m_{1}^E(\ch,\ch)\leq m_{1h}^E(\ch,\ch) \leq \kappa_1^* m_{1}^E(\ch,\ch),
\quad \forall \, \ch \in \M_h(E),\\
&\kappa_{2*} a_{0}^E(\vh,\vh)\leq a_{0h}^E(\vh,\vh) \leq \kappa_2^* a_{0}^E(\vh,\vh), \quad \forall\,\vh \in \V_h(E),\\
&\kappa_{3*} a_{1}^E(\ch,\ch)\leq a_{1h}^E(\ch,\ch) \leq \kappa_3^* \sigma^{-1}\|\ch\|_{1,E}^2,
\quad \forall \, \ch \in \M_h(E).
\end{align*}
\end{itemize}
\end{lemma}
\hspace*{-0.5cm}The discrete global bilinear and trilinear forms and load terms are defined as, for all $\uh,\wh,\vh \in \V_h$, $\bh,\dh,\ch \in \M_h$, 
\begin{align}
& m_{0h}(\uh,\vh) = \sum_{E \in\Th} m_{0h}^E(\uh,\vh),\quad m_{1h}(\bh,\ch)= \sum_{E \in\Th} m_{1h}^E(\bh,\ch),  \\
& a_{0h}(\uh,\vh) = \sum_{E \in\Th} a_{0h}^E(\uh,\vh),\quad a_{1h}(\bh,\ch)= \sum_{E \in\Th} a_{1h}^E(\bh,\ch),  \\
&c_{0h}(\uh,\wh,\vh) = \sum_{E\in\Th}c_{0h}^E(\uh,\wh,\vh) = \sum_{E\in\Th}(\Pi_{k-1}^{0,E}(\nabla \vh)\Pi_k^{0,E}\uh,\Pi_k^{0,E}\wh)_E,  \\
&c_{1h}(\bh,\dh,\vh)= \sum_{E\in\Th}c_{1h}^E(\bh,\dh,\vh) = \mu\sum_{E\in\Th}(\Pi_{k-1}^{0,E}\curl \bh \times \Pi_k^{0,E}\dh,\Pi_k^{0,E}\vh)_E,\\
&(\f_h,\vh) = \sum_{E\in\Th}(\Pi_k^{0,E}\f,\vh)_E, \quad (\g_h,\ch)= \sum_{E\in\Th}(\Pi_k^{0,E}\g,\ch)_E .
\end{align} 

Here, we introduce three useful lemmas for error analysis, which describe the interpolation error estimates and projection error estimate. The proof of Lemma \ref{interpolation-DFCVEM} can be found in \cite{interpolation-Stokes}. Lemma \ref{interpolation-CVEM} is from reference \cite{VEM-equivalent}. Lemma \ref{pi-estimate} is a classical result (see \cite{FEM2008}).
\begin{lemma}\label{interpolation-DFCVEM}
Let $\u \in \H^{k+1}(\O)$, $0< s \leq k$ and $\u_I \in \V_h$ be its degrees of freedom interpolation. Then, it holds
\begin{align*} 
\|\u - \u_I\|_{0,E}  + h|\u - \u_I|_{1,E} \leq C h^{s+1}|\v|_{s+1,E},\quad \forall E \in \Th .
\end{align*}
\end{lemma}
\begin{lemma}\label{interpolation-CVEM}
Let $\b \in \H^{k+1}(\O)$, $1 \leq s \leq k $ and $\b_I\in\M_h$ be its degrees of freedom interpolation. Then, it holds
\begin{align*}
\|\b - \b_I\|_{0,E} + h|\b - \b_I|_{1,E} \leq C h^{s+1}|\b|_{s+1,E}, \quad \forall \, E\in \Th.
\end{align*}
\end{lemma}
\begin{lemma} \label{pi-estimate}
Let $\v\in \H^{k+1}(E)$ and $0\leq s \leq k$, there exists a polynomial $\v_{\pi}\in [\P_k(E)]^2$ such that 
\begin{align*} 
\|\v - \v_{\pi}\|_{L^l(E)} + h |\v - \v_{\pi}|_{W^{1,l}(E)} \leq Ch^{s+1}|\v|_{W^{s+1,l}(E)}, \quad \forall \, 1\leq l \leq \infty, E\in\Th.
\end{align*}
\end{lemma}

\setcounter{equation}{0}
\section{Fully discrete BDF2 IMEX SAV scheme}
The MHD equations \refe{MHDeq} can be rewritten into the following equivalent system:
\begin{subequations}
\begin{align}
\u_t - \nu \Delta \u  + \nabla p + \frac{q(t)}{q(t)}\big((\nabla \u) \u -\mu \curl \b \times \b \big) = \f &,   \label{eqMHDeq1}\\
\mu\b_t + \sigma^{-1}\curl\curl\b - \frac{q(t)}{q(t)}\mu \curl(\u\times\b) = \g &, \label{eqMHDeq2}\\
\div \u = 0, \quad \div \b = 0 &, \label{eqMHDeq3} \\
\frac{dq(t)}{dt} = -\frac{1}{T}q(t) - \frac{1}{q(t)}\big(c_0(\u,\u,\u) +c_1(\b,\b,\u) - c_1(\b,\b,\u)\big)&,   \label{eqMHDeq4}
\end{align}
\end{subequations}
with the scalar auxiliary variable $q(t) = \exp(-\frac{t}{T})$. 

Let $\dta t = T/(N+1), t^n = n\dta t, 0\leq n \leq N+1,\u^n = \u(t^n),\b^n = \b(t^n),q^n = q(t^n)$. Due to BDF2 being a three-level scheme in time, let $(\uh^0,\bh^0,\qh^0) = (\u_I^0,\b_I^0,q^0)$, where $\u_I^0,\b_I^0$ are the interpolations of $\u^0$ and $\b^0$, and $(\uh^1,\bh^1,p_h^1,\qh^1)$ be the solution of the following BDF1 IMEX SAV scheme
\begin{subequations}
\begin{align}
m_{0h}(\frac{\uh^1 - \uh^0}{\dta t},\vh ) + a_{0h}(\frac{\uh^1 + \uh^0}{2},\vh) + d(\vh,p_h^1) - \frac{q_h^1}{q^1}\Big( c_{0h}(\uh^0,\uh^0,\vh)& \notag\\
  + c_{1h}(\bh^0,\bh^0,\vh)\Big) = (\fh^1,\vh),& \label{BDF1fdisMHDeq1}\\
m_{1h}(\frac{\bh^1 - \bh^0}{\dta t},\ch) + a_{1h}(\frac{\bh^1 + \bh^0}{2},\ch) + \frac{\qh^1}{q^1}c_{1h}(\ch,\bh^0,\uh^0) = (\gh^1,\ch),  & \label{BDF1fdisMHDeq2}\\
d(\uh^1,r_h) = 0, &\label{BDF1fdisMHDeq3}\\
\frac{\qh^1 - q_h^0}{\dta t} = -\frac{1}{T}q_h^1 - \frac{1}{q^1} \Big(c_{0h}(\uh^0,\uh^0,\uh^1) + c_{1h}(\bh^0,\bh^0,\uh^1)
 - c_{1h}(\bh^1,\bh^0,\uh^0) \Big)&,  \label{BDF1fdisMHDeq4} 
\end{align}
\end{subequations} 
for all $\vh \in \Vh$, $\ch \in \Mh$, $r_h \in Q_h$. We describe how to solve \refe{BDF1fdisMHDeq1}-\refe{BDF1fdisMHDeq4} efficiently. Set 
\begin{align}
\begin{aligned}
\uh^1 &= \u_{1,h}^1 + \xi^1 \u_{2,h}^1, \\
\bh^1 &= \b_{1,h}^1 + \xi^1 \b_{2,h}^1, \\
 p_h^1 &=  p_{1,h}^1 + \xi^1  p_{2,h}^1, \\
\end{aligned} \label{decoupset}
\end{align}
where $\xi^{1} = q_h^{1}/q^{1} $. By substituting \refe{decoupset} into \refe{BDF1fdisMHDeq1}-\refe{BDF1fdisMHDeq4}, we derive
\begin{align}
\begin{aligned}
m_{0h}(\frac{\u_{1,h}^1 - \uh^0}{\dta t},\vh) + a_{0h}(\frac{\u_{1,h}^0 + \uh^0}{2},\vh) + d(\vh,p_{1,h}^1) &= (\fh^1,\vh),\\
d(\u_{1,h}^1,\rh) &= 0,\\
m_{1h}(\frac{\b_{1,h}^1 - \bh^0}{\dta t},\ch) + a_{1h}(\frac{\b_{1,h}^1 + \bh^0}{2},\ch)  &= (\gh^1,\ch), 
\end{aligned} \label{decoupeq1} 
\end{align}
and
\begin{align}
\begin{aligned}
m_{0h}(\frac{\u_{2,h}^1}{\dta t},\vh) + a_{0h}(\frac{\u_{2,h}^1}{2},\vh) + d(\vh,p_{2,h}^1)- c_{0h}(\uh^0,\uh^0,\vh) - c_{1h}(\bh^0,\bh^0,\vh) &= 0, \\
d(\u_{2,h}^1,r_h) &=0, \\
m_{1h}(\frac{\b_{2,h}^1}{\dta t},\ch) + a_{1h}(\frac{\b_{2,h}^1}{2},\ch) + c_{1h}(\ch,\bh^0,\uh^0) &= 0.
\end{aligned} \label{decoupeq2}
\end{align}
By \refe{decoupeq1} and \refe{decoupeq2}, $\u_{i,h}^{1},p_{i,h}^{1}$ and $\b_{i,h}^{1}\,(i=1,2)$ can be obtained. Using \refe{decoupset}, we can rewrite \refe{BDF1fdisMHDeq4} as 
\begin{align}
\begin{aligned}
&\Big(\frac{\dta t + T}{T \dta t}  + \frac{1}{(q^1)^2}\Big(c_{0h}(\uh^0,\uh^0,\u_{2,h}^1) + c_{1h}(\bh^0,\bh^0,\u_{2,h}^1) - c_{1h}(\b_{2,h}^1,\bh^0,\uh^0) \Big)\Big)q_h^1 \\
&= \frac{\qh^0}{\dta t} - \frac{1}{q^1}\Big(c_{0h}(\uh^0,\uh^0,\u_{1,h}^1) + c_{1h}(\bh^0,\bh^0,\u_{1,h}^1) - c_{1h}(\b_{1,h}^1,\bh^0,\uh^0) \Big).
\end{aligned}  \label{decoupeq3}
\end{align}
Next, we prove that $\qh^1$ is uniquely solvable. Taking $\vh = \u_{2,h}^1, r_h = -p_{2,h}^1, \ch = \b_{2,h}^1$ in \refe{decoupeq2} and summing them up yield
\begin{align*}
&c_{0h}(\uh^0,\uh^0,\u_{2,h}^1) + c_{1h}(\bh^0,\bh^0,\u_{2,h}^1) - c_{1h}(\b_{2,h}^1,\bh^0,\uh^0)\\
& = \frac{1}{\dta t}m_{0h}(\u_{2,h}^1,\u_{2,h}^1) + \frac{1}{\dta t}m_{1h}(\b_{2,h}^1,\b_{2,h}^1) + \frac{1}{2} a_{0h}(\u_{2,h}^1,\u_{2,h}^1) +  \frac{1}{2} a_{1h}(\b_{2,h}^1,\b_{2,h}^1) \geq 0.
\end{align*}
Thus, $\qh^1$ and $\xi^1$ are uniquely solvable. Then we obtain $\uh^{1},p_h^{1}$ and $\bh^{1}$ by \refe{decoupset}. 

The fully discrete BDF2 IMEX SAV scheme can be given as follows: find $(\uh^{n+1},\bh^{n+1}$, $p_h^{n+1}$, $\qh^{n+1})$, $n=1, 2, \cdots,N$, such that
\begin{subequations}
\begin{align}
m_{0h}(\frac{3\uh^{n+1} - 4\uh^n + \uh^{n-1}}{2\dta t},\vh ) + a_{0h}(\uh^{n+1},\vh) + d(\vh,p_h^{n+1})  - \frac{\qh^{n+1}}{q^{n+1}}\Big( c_{0h}(\huh^n,\huh^n,\vh)
  \notag\\
  + c_{1h}(\hbh^n,\hbh^n,\vh)\Big) = (\fh^{n+1},\vh), \label{fdisMHDeq1}\\
m_{1h}(\frac{3\bh^{n+1} - 4\bh^n + \bh^{n-1}}{2\dta t},\ch) + a_{1h}(\bh^{n+1},\ch) + \frac{\qh^{n+1}}{q^{n+1}}c_{1h}(\ch,\hbh^n,\huh^n) = (\gh^{n+1},\ch),   \label{fdisMHDeq2}\\
d(\uh^{n+1},r_h) = 0, \label{fdisMHDeq3}\\
\frac{3\qh^{n+1} - 4\qh^n + \qh^{n-1}}{2\dta t} = -\frac{1}{T}\qh^{n+1} - \frac{1}{q^{n+1}} \Big(c_{0h}(\huh^n,\huh^n,\uh^{n+1}) + c_{1h}(\hbh^n,\hbh^n,\uh^{n+1})&\notag \\
 - c_{1h}(\bh^{n+1},\hbh^n,\huh^n) \Big),  \label{fdisMHDeq4} 
\end{align}
\end{subequations}
for all $\vh \in \Vh$, $\ch \in \Mh$, $r_h \in Q_h$, where $\overline{\v}^n = 2\v^n - \v^{n-1}$. Similar to the decoupling process of the BDF1 IMEX SAV scheme, we deduce $\qh^{n+1},\uh^{n+1}, p_h^{n+1}$ and $\bh^{n+1}$. 

\begin{lemma} \label{stabilityestimate}
Let $\f,\g \in L^2(0,T;\L^2(\O))$, $\forall\, n=1,...,N$, the following estimates hold
\begin{align}
\begin{aligned}
&\frac{1}{2}\Big( m_{0h}(\uh^1,\uh^1) + m_{1h}(\bh^1,\bh^1) + |\qh^1|^2 \Big) + \frac{1}{4} \dta t\Big(a_{0h}(\uh^1,\uh^1) + a_{1h}(\bh^1,\bh^1)  + \frac{4}{T} |\qh^1|^2 \Big)\\
& \leq \frac{1}{2}\Big(m_{0h}(\uh^0,\uh^0) + m_{1h}(\bh^0,\bh^0) + |\qh^0|^2\Big) +\frac{1}{2}\dta t\Big(a_{0h}(\uh^0,\uh^0) + a_{1h}(\bh^0,\bh^0) \Big) \\
&\quad + C\dta t\big( \|\f^1\|_0^2 + \|\g^1\|_0^2\big)
\end{aligned} \label{stabilityeq1}
\end{align}
and
\begin{align}
\begin{aligned}
&\frac{1}{2}\Big( m_{0h}(\uh^{n+1},\uh^{n+1}) + m_{0h}(\huh^{n+1},\huh^{n+1})+ m_{1h}(\bh^{n+1},\bh^{n+1}) + m_{1h}(\hbh^{n+1},\hbh^{n+1}) + |\qh^{n+1}|^2 \\
& + |\hqh^{n+1}|^2 \Big)  + \dta t \sum_{j=2}^{n+1}\Big(a_{0h}(\uh^j,\uh^j) + a_{1h}(\bh^j,\bh^j) + \frac{2}{T}|\qh^j|^2\Big)\\
& \leq  C\Big( m_{0h}(\uh^0,\uh^0) + m_{1h}(\bh^0,\bh^0)  + |\qh^0|^2 \Big) + C\dta t\Big(a_{0h}(\uh^0,\uh^0) + a_{1h}(\bh^0,\bh^0) \Big)  \\
&\quad + C\dta t\sum_{j=1}^{n+1}\big( \|\f^j\|_0^2 + \|\g^j\|_0^2\big) .
\end{aligned} \label{stabilityeq2}
\end{align}
\end{lemma}
\begin{proof}
Taking $(\vh,\ch,r_h) = \dta t(\uh^1,\bh^1,-p_h^1)$ in \refe{BDF1fdisMHDeq1}-\refe{BDF1fdisMHDeq3}, multiplying \refe{BDF1fdisMHDeq4} by $\dta t \qh^1$ and applying the identities 
\begin{align}
\begin{aligned}
2a(a-b) = a^2 - b^2 + (a-b)^2, \qquad  2a(a+b) = a^2 - b^2 + (a+b)^2,
\end{aligned} \label{eqeq1}
\end{align}
then, summing them up, using properties of $L^2$-projection and the stability in Lemma \ref{stalility-consistency}, we deduce
\begin{align*}
&\frac{1}{2}\Big( m_{0h}(\uh^1,\uh^1) + m_{1h}(\bh^1,\bh^1) + |\qh^{1}|^2 + m_{0h}(\uh^1 - \uh^0,\uh^1 - \uh^0) + m_{1h}(\bh^1 - \bh^0,\bh^1 - \bh^0) \\
&+  |\qh^1 - \qh^0|^2\Big) + \frac{1}{2}\dta t \Big(a_{0h}(\uh^1,\uh^1) + a_{1h}(\bh^1,\bh^1)  + \frac{2}{T} |\qh^1|^2 + a_{0h}(\uh^1 + \uh^0,\uh^1 + \uh^0) \\
& + a_{1h}(\bh^1+\bh^0,\bh^1+\bh^0)\Big)\\
&= \frac{1}{2}\big(m_{0h}(\uh^0,\uh^0) + m_{1h}(\bh^0,\bh^0) + |\qh^0|^2\big) + \frac{1}{2}\dta t\big(a_{0h}(\uh^0,\uh^0) + a_{1h}(\bh^0,\bh^0) \big)+ \dta t (\fh^1,\uh^1) \\
&\quad + \dta t(\gh^1,\bh^1 )  \\
& \leq \frac{1}{2}\big(m_{0h}(\uh^0,\uh^0) + m_{1h}(\bh^0,\bh^0) + |\qh^0|^2\big) +\frac{1}{2}\dta t\big(a_{0h}(\uh^0,\uh^0) + a_{1h}(\bh^0,\bh^0) \big) \\
&\quad + \frac{1}{4}\dta t \big(a_{0h}(\uh^1,\uh^1)+ a_{1h}(\bh^1,\bh^1) \big) + C\dta t\big( \|\f^1\|_0^2 + \|\g^1\|_0^2\big),
\end{align*}
which implies that the estimate \refe{stabilityeq1} holds. 

Next, we prove the estimate \refe{stabilityeq2}. Taking $(\vh,\ch,r_h) = 2\dta t(\uh^{n+1}, \bh^{n+1},-p_h^{n+1})$ in \refe{fdisMHDeq1}-\refe{fdisMHDeq3}, multiplying \refe{fdisMHDeq4} by $2\dta t \qh^{n+1}$ and applying the identity 
\begin{align}
2a(3a - 4b + c) = a^2 + (2a - b)^2 - b^2 - (2b - c)^2 + (a-2b + c)^2, \label{BDF2eq1}
\end{align}
then combining them and summing up from $j=2$ to $n+1$, we derive
\begin{align}
&\frac{1}{2}\Big( m_{0h}(\uh^{n+1},\uh^{n+1}) + m_{0h}(\huh^{n+1},\huh^{n+1})+ m_{1h}(\bh^{n+1},\bh^{n+1}) + m_{1h}(\hbh^{n+1},\hbh^{n+1}) + |\qh^{n+1}|^2 \notag \\
& + |\hqh^{n+1}|^2 \Big) + \frac{1}{2}\sum_{j=2}^{n+1}\Big( m_{0h}(\uh^j - \huh^{j-1},\uh^j - \huh^{j-1}) + m_{1h}(\bh^j - \hbh^{j-1},\bh^j - \hbh^{j-1})  \notag\\
&+ |\qh^j - \hqh^{j-1}|^2\Big) + 2\dta t \sum_{j=2}^{n+1}\Big(a_{0h}(\uh^j,\uh^j) + a_{1h}(\bh^j,\bh^j) + \frac{1}{T}|\qh^j|^2\Big)   \notag\\
&= \frac{1}{2}\Big( m_{0h}(\uh^1,\uh^1) + m_{0h}(\huh^1,\huh^1)+ m_{1h}(\bh^1,\bh^1) + m_{1h}(\hbh^1,\hbh^1) + |\qh^1|^2 + |\hqh^1|^2 \Big)  \notag\\
&\quad + 2\dta t\sum_{j=2}^{n+1} (\fh^j,\uh^j) +2\dta t\sum_{j=2}^{n+1} (\gh^j,\bh^j ) \notag \\
& \leq  C\Big( m_{0h}(\uh^0,\uh^0) + m_{1h}(\bh^0,\bh^0)  + |\qh^0|^2 \Big) + C\dta t\Big(a_{0h}(\uh^0,\uh^0) + a_{1h}(\bh^0,\bh^0) \Big)  \notag \\
&\quad + C\dta t\sum_{j=1}^{n+1}\big( \|\f^j\|_0^2 + \|\g^j\|_0^2\big)  + \dta t\sum_{j=2}^{n+1} \Big(a_{0h}(\uh^j,\uh^j) + a_{1h}(\bh^j,\bh^j) \Big).   \label{BDF2eq2}
\end{align}
In fact, by applying the properties of $L^2$-projection, Young's inequality, the stability from Lemma \ref{stalility-consistency}, and \refe{stabilityeq1}, the last inequality of \refe{BDF2eq2} is true. Consequently, we can infer \refe{stabilityeq2}. The proof is completed.
\end{proof}

\begin{lemma} \label{nc-bound} 
 Suppose $\u^0 \in \H^2(\O)\cap\V $ and  $\b^0 \in \H^2(\O)\cap\M $, then the solutions $(\uh^n,\bh^n,p_h^n,\qh^n)$ of \refe{BDF1fdisMHDeq1}-\refe{BDF1fdisMHDeq4} and \refe{fdisMHDeq1}-\refe{fdisMHDeq4}, and $(\uh^0,\bh^0,q_h^0)$ satisfy
\begin{align*}
&\|\uh^n\|_0^2 + \|\bh^n\|_0^2  + |q_h^n|^2 \leq C, \quad \forall\, 0\leq n \leq N+1,\\
\|\nabla \uh^0\|_0^2 + \|\bh^0\|_1^2 &\leq C ,\quad \dta t \sum_{j=1}^{n}\big(\|\nabla \uh^j\|_0^2 + \|\bh^j\|_1^2\big)  \leq C, \quad \forall \,1 \leq n \leq N+1  .
\end{align*}
\end{lemma}
\begin{proof}
From Lemmas \ref{interpolation-DFCVEM} and \ref{interpolation-CVEM}, we observe that the above inequalities are true for $n=0$. Utilizing Lemma \ref{stabilityestimate} and the stablity in Lemma \ref{stalility-consistency}, for $1\leq n \leq N+1$, we can obtain the desire results.
\end{proof}
\begin{remark} \label{re-uh0}
Note that, by Lemmas \ref{interpolation-DFCVEM}-\ref{pi-estimate}, it is not difficult to find that under assumption of Lemma \ref{nc-bound}, $\|\Pi_k^{0,E}\uh^0\|_{L^{\infty}(E)}$ and $\|\Pi_k^{0,E}\bh^0\|_{L^{\infty}(E)}$ are bounded for all $E \in \Th$.
\end{remark}

\begin{lemma} \label{nolinear_estimate}
Let $\r,\w \in \H^{k+1}(\O)\cap \V, \a,\d \in \H^{k+1}(\O)\cap\M,\v \in \V,\c \in \M,$ and $\rrh,\wh \in \Vh,\ah,\dh \in \Mh$. Let $\alpha = \max_{E\in\Th}\|\Pi_k^{0,E}\rrh\|_{L^{\infty}(E)} $ and $\beta = \max_{E\in\Th} \|\Pi_k^{0,E}\dh\|_{L^{\infty}(E)} $, then the following estimates hold
\begin{subequations}
\begin{align}
\left| c_0(\r,\w,\v) - c_{0h}(\r,\w,\v) \right|  &\leq  Ch^k \|\r\|_{k+1}\|\w\|_{k+1}\|\nabla \v\|_0,  \label{nolineareq1}\\
\left|c_1(\a,\d,\v) - c_{1h}(\a,\d,\v)\right|    &\leq  Ch^k \|\a\|_{k+1}\|\d\|_k\|\nabla \v\|_0,   \label{nolineareq2}\\
\left|c_1(\c,\d,\r) - c_{1h}(\c,\d,\r)\right|     &\leq  Ch^k \|\d\|_{k+1}\|\r\|_{k+1}\|\c\|_1,  \label{nolineareq3} \\
\left|c_{0h}(\r,\w,\v) - c_{0h}(\rrh,\wh,\v)\right| &\leq C\big(\|\r - \rrh\|_0\|\w\|_{L^{\infty}} + \alpha\|\w - \wh\|_0 \big)\|\nabla \v\|_0,  \label{nolineareq4} \\
\left|c_{1h}(\a,\d,\v) - c_{1h}(\ah,\dh,\v)\right| &\leq  C \|\a\|_2\|\d - \dh\|_0\|\nabla \v\|_0 + C\beta\|\a - \ah\|_1 \| \v\|_0, \label{nolineareq5}    \\
\left|c_{1h}(\c,\d,\r) - c_{1h}(\c,\dh,\rrh)\right| &\leq C\big(\|\d -\dh\|_0\|\r\|_{L^{\infty}} + \beta\|\r -\rrh\|_0 \big)\|\c\|_1.  \label{nolineareq6}
\end{align}
\end{subequations} 
\end{lemma}
\begin{proof}
For the sake of simplicity, we only provide the proof of \refe{nolineareq1} and \refe{nolineareq4}. Applying the properties of $L^2$-projection and Sobolev embedding Theorem yields
\begin{align*}
&\left| c_0(\r,\w,\v) - c_{0h}(\r,\w,\v) \right|\\
=& \Big|\sum_{E \in \Th} \Big(((\nabla\v)\r,\w)_E - (\Pi_{k-1}^{0,E}(\nabla\v)\Pi_k^{0,E}\r,\Pi_k^{0,E}\w)_E \Big) \Big| \\
=& \Big|\sum_{E \in \Th} \Big((\I - \Pi_{k-1}^{0,E})(\nabla \v),(\I - \Pi_{k-1}^{0,E})(\w \otimes \r^{\mathrm{T}}))_E + (\Pi_{k-1}^{0,E}(\nabla \v) (\I -\Pi_k^{0,E})\r,\w)_E \\
& \quad + (\Pi_{k-1}^{0,E}(\nabla \v)\Pi_k^{0,E}\r,(\I - \Pi_k^{0,E})\w)_E\Big) \Big| \\
 \leq& Ch^k \|\w \otimes \r^{\mathrm{T}}\|_k \|\nabla \v\|_0 + Ch^k \|\r\|_{k+1}\|\w\|_1\|\nabla \v\|_0 + Ch^k \|\r\|_1\|\w\|_{k+1}\|\nabla \v\|_0 \\
 \leq& Ch^k \|\r\|_{k+1}\|\w\|_{k+1}\|\nabla \v\|_0.
\end{align*}
According to H$\ddot{\mathrm{o}}$lder inequality and the properties of $L^2$-projection, we observe that
\begin{align*}
&\left| c_{0h}(\r,\w,\v) - c_{0h}(\rrh,\wh,\v) \right| \\
=& \left| c_{0h}(\r - \rrh,\w,\v) - c_{0h}(\rrh,\w - \wh,\v)  \right| \\
\leq & C\|\r - \rrh\|_0\|\w\|_{L^{\infty}}\|\nabla \v\|_0 + \sum_{E \in \Th}\|\Pi_k^{0,E}\rrh\|_{L^{\infty}(E)}\|\w -\wh\|_{0,E} \|\nabla \v\|_{0,E} \\
\leq & C\big(\|\r - \rrh\|_0\|\w\|_{L^{\infty}} + \alpha\|\w - \wh\|_0 \big)\|\nabla \v\|_0.
\end{align*}
With just a few modifications, the rest of Lemma \ref{nolinear_estimate} can be got. The proof is completed.
\end{proof}

We end this section by introducing the Gronwall's inequalities from reference \cite{FEM2016}.
\begin{lemma} \label{gronwelleq}
Let $\varphi_n,\theta_n,\phi_n,\gamma_n, \dta t, C_0$ be non-negative numbers for $n \geq 1$ and the inequality 
\begin{align*}
\varphi_{N+1} + \dta t \sum_{n=1}^{N+1} \theta_n \le \dta t \sum_{n=1}^{N+1} \gamma_n \varphi_n + \dta t \sum_{n=1}^{N+1} \phi_n + C_0, \quad \forall\, N \geq 0
\end{align*}
hold. If $\dta t \gamma_n \leq 1$ for all $n=1,\cdots,N+1$, then
\begin{align*}
\varphi_{N+1} + \dta t \sum_{n=1}^{N+1} \theta_n \le \exp\Big(\dta t \sum_{n=1}^{N+1} \frac{\gamma_n}{1 - \dta t \gamma_n} \Big)\Big(\dta t \sum_{n=1}^{N+1} \phi_n + C_0\Big), \quad \forall\, N \geq 0.
\end{align*}

If the inequality
\begin{align*}
\varphi_{N+1} + \dta t \sum_{n=1}^{N+1} \theta_n \le \dta t \sum_{n=1}^N \gamma_n \varphi_n + \dta t \sum_{n=1}^{N+1} \phi_n + C_0, \quad \forall\, N \geq 1
\end{align*} 
is given, then it holds
\begin{align*}
\varphi_{N+1} + \dta t \sum_{n=1}^{N+1} \theta_n \le \exp\Big(\dta t \sum_{n=1}^N \gamma_n \Big)\Big(\dta t \sum_{n=1}^{N+1} \phi_n + C_0\Big), \quad \forall\, N \geq 1.
\end{align*}
\end{lemma}

\setcounter{equation}{0}
\section{Error analysis}
This section will establish error estimates. First of all, we provide the one between the exact solution $(\u,\b,p,q)$ of \refe{eqMHDeq1}-\refe{eqMHDeq4} and the numerical solution $(\uh^1,\bh^1,p_h^1,\qh^1)$ of \refe{BDF1fdisMHDeq1}-\refe{BDF1fdisMHDeq4}. Then the error estimate between the exact solution $(\u,\b,p,q)$ to \refe{eqMHDeq1}-\refe{eqMHDeq4} and the numerical solution $(\uh^{n+1},\bh^{n+1},p_h^{n+1},\qh^{n+1})$ to \refe{fdisMHDeq1}-\refe{fdisMHDeq4} is given. In order to derive the error estimates, we make the following regular assumption of solution to problem \refe{weakMHDeq}. 
\begin{assumption} \label{regassume}
The exact solution of \refe{weakMHDeq} satisfies
\begin{align*}
\begin{aligned}
&\u_{ttt},\b_{ttt} \in L^2(0,T;\L^2(\O)), \quad \u_{tt},\b_{tt} \in L^2(0,T;\H^1(\O)),\quad \u_t,\b_t \in L^{2}(0,T;\H^{k+1}(\O)),\\
&\u,\b \in L^{\infty}(0,T;\H^{k+1}(\O)) ,\quad  \u,\b \in L^{\infty}(0,T;\L^{\infty}(\O)),\quad p \in L^2(0,T;H^k(\O)). 
\end{aligned} 
\end{align*}
\end{assumption}
For the sake of convenience, for all $0\leq n \leq N+1$, we set
\begin{align*}
\u^n - \uh^n &= (\u^n - \u_I^n) + (\u_I^n - \uh^n) = \rho_{\u}^n + e_{\u}^n, \\
\b^n - \bh^n &= (\b^n - \b_I^n) + (\b_I^n - \bh^n) = \rho_{\b}^n + e_{\b}^n, \quad e_q^n = q^n - \qh^n,
\end{align*}
where $\u_I^n \in \Zh ,\b_I^n\in \Mh$ are the interpolations of $\u^n$ and $\b^n$. Let $\u_{\pi}^{n},\b_{\pi}^{n}$, given in Lemma \ref{pi-estimate}, be approximations of $\u^n$ and $\b^n$. From \refe{BDF1fdisMHDeq3}, \refe{fdisMHDeq3} and the fact $\div \V_h \in Q_h$, we see that $\div \uh^n = 0$ and $\div e_{\u}^n = 0$.

Taking the average of variational form of \refe{eqMHDeq1}-\refe{eqMHDeq2} at $t^0$ and $t^1$, choosing the mean of \refe{eqMHDeq4} at $t^0$ and $t^1$, we have
\begin{align}
\begin{aligned}
m_0(\frac{\u_t^0 + \u_t^1}{2},\vh) + a_0(\frac{\u^0 + \u^1}{2},\vh ) +  d(\vh,\frac{p^0 + p^1}{2}) - \frac{1}{2}\big(c_0(\u^0,\u^0,\vh) + c_0(\u^1,\u^1,\vh) \\
 + c_1(\b^0,\b^0,\vh) + c_1(\b^1,\b^1,\vh)\big) = \frac{1}{2}(\f^0 + \f^1,\vh),\\
m_1(\frac{\b_t^0 + \b_t^1}{2},\ch) + a_1(\frac{\b^0 + \b^1}{2},\ch) + \frac{1}{2}\big( c_1(\ch,\b^0,\u^0) + c_1(\ch,\b^1,\u^1)\big) = \frac{1}{2}(\g^0 + \g^1,\ch),\\
\frac{q_t^0 + q_t^1}{2}= - \frac{q^0 + q^1}{2T}- \frac{1}{2q^0}\big(c_0(\u^0,\u^0,\u^0) + c_1(\b^0,\b^0,\u^0) - c_1(\b^0,\b^0,\u^0) \big) \\
-  \frac{1}{2q^1}\big(c_0(\u^1,\u^1,\u^1) + c_1(\b^1,\b^1,\u^1) - c_1(\b^1,\b^1,\u^1) \big).
\end{aligned} \label{1ceq1}
\end{align}
Subtrating \refe{BDF1fdisMHDeq1}-\refe{BDF1fdisMHDeq2} and \refe{BDF1fdisMHDeq4} from \refe{1ceq1}, then employing the identity
\begin{align*}
\frac{\v_t^1 + \v_t^0}{2} - \frac{\v^1 - \v^0}{\dta t} = \frac{1}{2\dta t}\int_{t^0}^{t^1}(t^1 - t)(t - t^0)\v_{ttt}\,dt, \quad \forall \,\v_{ttt} \in L^2(t^0,t^1;\L^2)
\end{align*}
and the consistency in Lemma \ref{stalility-consistency}, we derive
\begin{align}
&m_0(\frac{(\u^1 -\u^0 ) - (\u_{\pi}^1 -\u_{\pi}^0 )}{\dta t},\vh) + m_{0h}(\frac{(\u_{\pi}^1 - \u_{\pi}^0) - (\u_I^1 - \u_I^0)}{\dta t},\vh) + m_{0h}(\frac{e_{\u}^1 - e_{\u}^0}{\dta t},\vh) \notag\\
& +\sum_{E \in \Th}a_0^E(\frac{(\u^1 + \u^0) - (\u_{\pi}^1 + \u_{\pi}^0)}{2} ,\vh) + a_{0h}(\frac{(\u_{\pi}^1 + \u_{\pi}^0) -(\u_I^1 + \u_I^0)}{2},\vh ) \notag\\
&+ a_{0h}(\frac{e_{\u}^1 + e_{\u}^0}{2},\vh) - \Big\{c_0(\u^0,\u^0,\vh) - c_{0h}(\u^0,\u^0,\vh) + c_{0h}(\u^0,\u^0,\vh) - c_{0h}(\uh^0,\uh^0,\vh)  \notag\\
& + \frac{e_q^1}{q^1}c_{0h}(\uh^0,\uh^0,\vh)\Big\}   - \Big\{c_1(\b^0,\b^0,\vh) - c_{1h}(\b^0,\b^0,\vh) + c_{1h}(\b^0,\b^0,\vh)  \notag\\
&-  c_{1h}(\bh^0,\bh^0,\vh) + \frac{e_q^1}{q^1}c_{1h}(\bh^0,\bh^0,\vh)\Big\} - (\f^1 - \fh^1,\vh) + d(\vh,\frac{p^0 + p^1}{2} - p_h^1) = (G_1,\vh), \label{1ceq2}\\
&m_1(\frac{(\b^1 -\b^0 ) - (\b_{\pi}^1 -\b_{\pi}^0 )}{\dta t},\ch) + m_{1h}(\frac{(\b_{\pi}^1 - \b_{\pi}^0) - (\b_I^1 - \b_I^0)}{\dta t},\ch) + m_{1h}(\frac{e_{\b}^1 - e_{\b}^0}{\dta t},\ch) \notag\\
& + \sum_{E\in \Th} a_1^E(\frac{(\b^1 + \b^0) - (\b_{\pi}^1 + \b_{\pi}^0)}{2},\ch) + a_{1h}(\frac{(\b_{\pi}^1 + \b_{\pi}^0) -(\b_I^1 + \b_I^0)}{2},\ch) \notag\\
&+ a_{1h}(\frac{e_{\b}^1 + e_{\b}^0}{2},\ch)  + \Big\{c_1(\ch,\b^0,\u^0) - c_{1h}(\ch,\b^0,\u^0) + c_{1h}(\ch,\b^0,\u^0) - c_{1h}(\ch,\bh^0,\uh^0) \notag\\
&+ \frac{e_q^1}{q^1}c_{1h}(\ch,\bh^0,\uh^0) \Big\}  - (\g^1 - \gh^1,\ch) = (G_2,\ch), \label{1ceq3} \\
& \frac{e_q^1 - e_q^0}{\dta t} + \frac{e_q^1}{T} + \frac{1}{q^1}\Big\{ c_0(\u^0,\u^0,\u^1) - c_{0h}(\u^0,\u^0,\u^1) + c_{0h}(\u^0,\u^0,\u^1)  - c_{0h}(\uh^0,\uh^0,\u^1) \notag\\
&+ c_{0h}(\uh^0,\uh^0,\rho_{\u}^1) + c_{0h}(\uh^0,\uh^0,e_{\u}^1)\Big\} + \frac{1}{q^1}\Big\{c_1(\b^0,\b^0,\u^1) - c_{1h}(\b^0,\b^0,\u^1)+ c_{1h}(\b^0,\b^0,\u^1) \notag\\
& - c_{1h}(\bh^0,\bh^0,\u^1) + c_{1h}(\bh^0,\bh^0,\rho_{\u}^1) + c_{1h}(\bh^0,\bh^0,e_{\u}^1) \Big\} - \frac{1}{q^1}\Big\{c_1(\b^1,\b^0,\u^0) - c_{1h}(\b^1,\b^0,\u^0) \notag\\
& + c_{1h}(\b^1,\b^0,\u^0) - c_{1h}(\b^1,\bh^0,\uh^0)  + c_{1h}(\rho_{\b}^1,\bh^0,\uh^0) + c_{1h}(e_{\b}^1,\bh^0,\uh^0) \Big\} = G_3,  \label{1ceq4}
\end{align}
where
\begin{align}
(G_1,\vh) 
&= - \frac{1}{2\dta t} \int_{t^0}^{t^1}(t^1 - t)(t- t^0)(\u_{ttt},\vh)\,dt + \frac{1}{2}\Big\{c_0(\int_{t^0}^{t^1}\u_t\,dt,\u^1,\vh) \notag\\
& \quad+ c_0(\u^0,\int_{t^0}^{t^1}\u_t\,dt,\vh) \Big\} 
 + \frac{1}{2}\Big\{c_1(\int_{t^0}^{t^1}\b_t\,dt,\b^1,\vh) + c_1(\b^0,\int_{t^0}^{t^1}\b_t\, dt,\vh) \Big\} \notag\\
 &\quad - \frac{1}{2}(\int_{t^0}^{t^1}\f_t,\vh)\,dt,  \label{1ceq5}\\
(G_2,\ch) 
&= -\frac{\mu}{2\dta t}\int_{t^0}^{t^1}(t^1 - t)(t- t^0)(\b_{ttt},\ch)\,dt - \frac{1}{2}\Big\{ c_1(\ch,\int_{t^0}^{t^1}\b_t\,dt,\u^1) \notag\\
& \quad+ c_1(\ch,\b^0,\int_{t^0}^{t^1}\u_t\,dt) \Big\} - \frac{1}{2}(\int_{t^0}^{t^1}\g_t,\ch)\,dt,  \label{1ceq6} \\
G_3 
&= -\frac{1}{2\dta t}\int_{t^0}^{t^1}(t^1 - t)(t - t^0)q_{ttt}\,dt + \frac{1}{2T}\int_{t^0}^{t^1} q_t\,dt - \frac{1}{2q^1}\Big\{ \frac{1}{q^0}c_0(\u^0,\u^0,\u^0)\int_{t^0}^{t^1}q_t \,dt  \notag\\
& \quad + c_0(\int_{t^0}^{t^1}\u_t\,dt,\u^1,\u^1) + c_0(\u^0,\int_{t^0}^{t^1} \u_t\,dt,\u^1) -  c_0(\u^0,\u^0,\int_{t^0}^{t^1} \u_t\,dt)\Big\} \notag\\
& \quad  - \frac{1}{2q^1}\Big\{ \frac{1}{q^0}c_1(\b^0,\b^0,\u^0)\int_{t^0}^{t^1} q_t\,dt + c_1(\int_{t^0}^{t^1}\b_t\,dt,\b^1,\u^1)   + c_1(\b^0,\int_{t^0}^{t^1}\b_t\,dt,\u^1) \notag\\
& \quad- c_1(\b^0,\b^0,\int_{t^0}^{t^1}\u_t\,dt)  \Big\}  + \frac{1}{2q^1}\Big\{\frac{1}{q^0}c_1(\b^0,\b^0,\u^0)\int_{t^0}^{t^1}q_t\,dt + c_1(\b^1,\b^1,\int_{t^0}^{t^1}\u_t\,dt) \notag\\
& \quad+ c_1(\b^1,\int_{t^0}^{t^1}\b_t\,dt,\u^0) - c_1(\int_{t^0}^{t^1}\b_t\,dt,\b^0,\u^0)\Big\}. \label{1ceq7}
\end{align}
\begin{lemma} \label{G_estimate}
Suppose that the Assumption \ref{regassume} is satisfied, the following estimate holds
\begin{align*}
\|G_1\|_0 + \|G_2\|_0 + |G_3| \leq C \dta t.
\end{align*} 
\end{lemma}
\begin{proof}
From \refe{1ceq5} and an integration by parts, we observe that
\begin{align*}
\|G_1\|_0 
&= \sup_{0 \neq \vh \in \L^2(\O)} \frac{(G_1,\vh)}{\|\vh\|_0} \\
& \leq C (\dta t)^{\frac{3}{2}}\Big(\int_{t^0}^{t^1} \|\u_{ttt}\|_0^2\,dt \Big)^{\frac{1}{2}} + C\dta t \|\u_t\|_{L^{\infty}(t^0,t^1;\H^1(\O))}\big(\|\u^0\|_{L^{\infty}} + \|\u^1\|_2 \big) \\
& \quad + C\dta t  \|\b_t\|_{L^{\infty}(t^0,t^1;\H^1(\O))}\big(\|\b^1\|_{L^{\infty}} + \|\b^0\|_2 \big) + C \dta t \|\f_t\|_{L^{\infty}(t^0,t^1;\L^2(\O))}.
\end{align*}
According to \refe{1ceq6}, an integration by parts and identity \refe{curl-eq1}, the following bound can be inferred
\begin{align*}
\|G_2\|_0 
&= \sup_{0 \neq \ch \in \L^2(\O)} \frac{(G_2,\ch)}{\|\ch\|_0} \\
&\leq C(\dta t)^{\frac{3}{2}}\Big(\int_{t^0}^{t^1}\|\b_{ttt}\|_0^2 \,dt \Big)^{\frac{1}{2}} + C\dta t  \|\b_t\|_{L^{\infty}(t^0,t^1;\H^1(\O))}\big(\|\u^1\|_{L^{\infty}} + \|\u^1\|_2\big) \\
& \quad + C\dta t  \|\u_t\|_{L^{\infty}(t^0,t^1;\H^1(\O))}\big(\|\b^0\|_2 + \|\b^0\|_{L^{\infty}} \big) + C\dta t \|\g_t\|_{L^{\infty}(t^0,t^1;\L^2(\O))}.
\end{align*}
Analogously,
\begin{align*}
|G_3| 
&\leq C(\dta t)^{\frac{3}{2}}\Big(\int_{t^0}^{t^1}|q_{ttt}|\,dt \Big)^{\frac{1}{2}} + C\dta t\max_{t^0 \leq t \leq t^1}|q_t|\big(1 + \|\u^0\|_1^3 + \|\b^0\|_1^2\|\u^0\|_1 \big) \\
& \quad + C\dta t  \|\u_t\|_{L^{\infty}(t^0,t^1;\H^1(\O))}\big(\|\u^1\|_1^2 + \|\u^0\|_1\|\u^1\|_1 + \|\u^0\|_1^2 + \|\b^0\|_1^2 + \|\b^1\|_1^2\big) \\
& \quad + C \dta t  \|\b_t\|_{L^{\infty}(t^0,t^1;\H^1(\O))}\big(\|\b^1\|_1\|\u^1\|_1 + \|\b^0\|_1\|\u^1\|_1 + \|\b^1\|_1\|\u^0\|_1 + \|\b^0\|_1\|\u^0\|_1 \big).
\end{align*}
Combining the above inequalities, the proof is completed.
\end{proof}
\begin{theorem} \label{1c-error-estimate} 
Assume that the Assumption \ref{regassume} is satisfied, then there holds 
\begin{align*}
\frac{1}{2}\big( \|e_{\u}^1\|_0^2 + \|e_{\b}^1\|_0^2 + |e_q^1|^2 \big) + \frac{1}{4}\dta t \big(\|\nabla e_{\u}^1\|_0^2 + \|e_{\b}^1\|_1^2 + \frac{2}{T}|e_q^1|^2 \big) \leq C \big((\dta t)^4  + h^{2k}\big).
\end{align*}
Further,
\begin{align*}
&\frac{1}{2}\big(\|\u^1 - \uh^1\|_0^2 + \|\b^1 - \bh^1\|_0^2 + |q^1 - \qh^1|^2 \big) + \frac{1}{4}\dta t \big(\|\nabla(\u^1 - \uh^1)\|_0^2 \\
& + \|\b^1 - \bh^1\|_1^2 + \frac{2}{T}|q^1 - \qh^1|^2 \big)  \leq C\big( (\dta t)^4 + h^{2k} \big).
\end{align*}
\end{theorem}
\begin{proof}
Taking $(\vh,\ch) = 2\dta t(e_{\u}^1,e_{\b}^1)$ in \refe{1ceq2}-\refe{1ceq3}, multiplying \refe{1ceq4} by $2\dta t e_q^1$ and using \refe{eqeq1} yield
\begin{align}
&m_{0h}(e_{\u}^1,e_{\u}^1) + m_{0h}(e_{\u}^1 - e_{\u}^0,e_{\u}^1 - e_{\u}^0) +  m_{1h}(e_{\b}^1,e_{\b}^1)  + m_{1h}(e_{\b}^1 - e_{\b}^0,e_{\b}^1 - e_{\b}^0) + |e_q^1|^2 \notag\\
&+ |e_q^1 - e_q^0|^2  + \frac{1}{2}\dta t \big(a_{0h}(e_{\u}^1,e_{\u}^1) + a_{0h}(e_{\u}^1 + e_{\u}^0,e_{\u}^1 + e_{\u}^0) + a_{1h}(e_{\b}^1,e_{\b}^1) \notag\\
& + a_{1h}(e_{\b}^1 + e_{\b}^0,e_{\b}^1 + e_{\b}^0) + \frac{4}{T}|e_q^1|^2 \big)  \notag \\
&= \Big\{m_{0h}(e_{\u}^0,e_{\u}^0) + m_{1h}(e_{\b}^0,e_{\b}^0) + |e_q^0|^2 + \frac{1}{2}\dta t \big(a_{0h}(e_{\u}^0,e_{\u}^0) + a_{1h}(e_{\b}^0,e_{\b}^0) \big) \notag\\
& \quad + m_0((\u_{\pi}^1 -\u_{\pi}^0) - (\u^1 -\u^0 ),2e_{\u}^1) + m_{0h}((\u_I^1 - \u_I^0) - (\u_{\pi}^1 - \u_{\pi}^0),2e_{\u}^1) \notag\\
& \quad + m_1((\b_{\pi}^1 -\b_{\pi}^0 ) - (\b^1 -\b^0 ),2e_{\b}^1) + m_{1h}((\b_I^1 - \b_I^0) - (\b_{\pi}^1 - \b_{\pi}^0),2e_{\b}^1)  \notag\\
& \quad  + \dta t\sum_{E \in \Th}a_0^E((\u_{\pi}^1 + \u_{\pi}^0) - (\u^1 + \u^0) ,e_{\u}^1) + \dta t a_{0h}((\u_I^1 + \u_I^0) - (\u_{\pi}^1 + \u_{\pi}^0),e_{\u}^1 )\notag\\
& \quad + \dta t\sum_{E\in \Th} a_1^E((\b_{\pi}^1 + \b_{\pi}^0) - (\b^1 + \b^0),e_{\b}^1) + \dta t a_{1h}((\b_I^1 + \b_I^0) - (\b_{\pi}^1 + \b_{\pi}^0),e_{\b}^1) \Big\} \notag\\ 
& \quad + \dta t \big\{c_0(\u^0,\u^0,2e_{\u}^1) - c_{0h}(\u^0,\u^0,2e_{\u}^1) + c_{0h}(\u^0,\u^0,2e_{\u}^1) - c_{0h}(\uh^0,\uh^0,2e_{\u}^1) \big\}  \notag\\
& \quad + \dta t \big\{c_1(\b^0,\b^0,2e_{\u}^1) - c_{1h}(\b^0,\b^0,2e_{\u}^1) + c_{1h}(\b^0,\b^0,2e_{\u}^1) - c_{1h}(\bh^0,\bh^0,2e_{\u}^1) \big\} \notag\\
& \quad - \dta t \big\{c_1(2e_{\b}^1,\b^0,\u^0) - c_{1h}(2e_{\b}^1,\b^0,\u^0) + c_{1h}(2e_{\b}^1,\b^0,\u^0) - c_{1h}(2e_{\b}^1,\bh^0,\uh^0)  \big\} \notag\\
&  \quad - 2\dta t\frac{e_q^1}{q^1}\Big\{ c_0(\u^0,\u^0,\u^1) - c_{0h}(\u^0,\u^0,\u^1) + c_{0h}(\u^0,\u^0,\u^1)  - c_{0h}(\uh^0,\uh^0,\u^1) \notag\\
& \quad + c_{0h}(\uh^0,\uh^0,\rho_{\u}^1)\Big\} -2\dta t \frac{e_q^1}{q^1}\Big\{c_1(\b^0,\b^0,\u^1) - c_{1h}(\b^0,\b^0,\u^1)+ c_{1h}(\b^0,\b^0,\u^1) \notag\\
& \quad - c_{1h}(\bh^0,\bh^0,\u^1) + c_{1h}(\bh^0,\bh^0,\rho_{\u}^1) \Big\} + 2\dta t \frac{e_q^1}{q^1}\Big\{c_1(\b^1,\b^0,\u^0) - c_{1h}(\b^1,\b^0,\u^0)  \notag\\
& \quad+ c_{1h}(\b^1,\b^0,\u^0) - c_{1h}(\b^1,\bh^0,\uh^0)  + c_{1h}(\rho_{\b}^1,\bh^0,\uh^0)  \Big\} + \dta t \big\{  (\f^1 - \fh^1,2e_{\u}^1) \notag\\
& \quad + (\g^1 - \gh^1,2e_{\b}^1) \big\} + \dta t \big\{(G_1,2e_{\u}^1) + (G_2,2e_{\b}^1) + 2e_q^1 \cdot G_3  \big\} \notag\\
&= \Gamma_1 + \Gamma_2 + \Gamma_3 + \Gamma_4 + \Gamma_5 + \Gamma_6 + \Gamma_7 + \Gamma_8 + \Gamma_9. \label{1c_estimateeq1}
\end{align}  
Here and after, we will always mark the components of long equations  according to the braces.  Next, we will give the estimates of $\Gamma_1$-$\Gamma_9$. Making use of Lemmas \ref{interpolation-DFCVEM}-\ref{pi-estimate}, the stability in Lemma \ref{stalility-consistency} and Cauchy–Schwarz inequality, we have
\begin{align}
\begin{aligned}
\Gamma_1 
&\leq Ch^{2k+2}\big(\|\u^0\|_{k+1}^2 + \|\b^0\|_{k+1}^2 \big) + C  h^{2k}\dta t\big(\|\u^0\|_{k+1}^2 + \|\b^0\|_{k+1}^2 \big) + Ch^{2k+2}\big(\|\u^0\|_{k+1}^2 \\
& \quad+ \|\u^1\|_{k+1}^2  + \|\b^0\|_{k+1}^2 + \|\b^1\|_{k+1}^2\big) + \frac{1}{4}m_{0h}(e_{\u}^1,e_{\u}^1) + \frac{1}{4}m_{1h}(e_{\b}^1,e_{\b}^1) + C  h^{2k}\dta t\big(\|\u^0\|_{k+1}^2 \\
& \quad +\|\u^1\|_{k+1}^2  + \|\b^0\|_{k+1}^2 +\|\b^1\|_{k+1}^2 \big) + \frac{1}{16}\dta t a_{0h}(e_{\u}^1,e_{\u}^1) + \frac{1}{12}\dta t a_{1h}(e_{\b}^1,e_{\b}^1).
\end{aligned} \label{1c_estimateeq2}
\end{align}
From \refe{nolineareq1}, \refe{nolineareq4} and Remark \ref{re-uh0}, we observe that
\begin{align}
\begin{aligned}
\Gamma_2 
&\leq  C\dta t \big(h^k\|\u^0\|_{k+1}^2 + \|\u^0 - \uh^0\|_0\|\u^0\|_{L^{\infty}} + \|\u^0 - \uh^0\|_0\big)\|\nabla e_{\u}^1\|_0 \\
&\leq C h^{2k}\dta t \big(\|\u^0\|_{k+1}^4 + h^2\|\u^0\|_{k+1}^2 \big) + \frac{1}{16}\dta t a_{0h}(e_{\u}^1,e_{\u}^1).
\end{aligned} \label{1c_estimateeq3}
\end{align}
Analogously, by Lemmas \ref{nc-bound}-\ref{nolinear_estimate}  and Remark \ref{re-uh0}, the following estimates can be derived
\begin{align}
\begin{aligned}
\Gamma_3 
&\leq Ch^{2k} \dta t\big(1 + h^2\|\b^0\|_2^2 + \|\b^0\|_k^2\big)\|\b^0\|_{k+1}^2 +\frac{1}{16}\dta t a_{0h}(e_{\u}^1,e_{\u}^1),  \\
\Gamma_4 
&\leq Ch^{2k}\dta t \big(h^2 + \|\b^0\|_{k+1}^2\big) \big(h^2 +\|\u^0\|_{k+1}^2\big) + \frac{1}{12}\dta ta_{1h}(e_{\b}^1,e_{\b}^1),\\ 
\Gamma_5 
& \leq C h^{2k}\dta t\big(\|\u^0\|_{k+1}^4 +  h^2\|\u^0\|_{k+1}^2\big)\|\nabla \u^1\|_0^2 + Ch^{2k}\dta t\|\u^1\|_{k+1}^2 +\frac{1}{2T}\dta t |e_q^1|^2,\\
\Gamma_6 
&\leq Ch^{2k} \dta t\big(1 + h^2\|\b^0\|_2^2 + \|\b^0\|_k^2\big)\|\b^0\|_{k+1}^2\|\nabla \u^1\|_0^2 + Ch^{2k+2}\dta t \|\u^1\|_{k+1}^2 +\frac{1}{2T}\dta t |e_q^1|^2, \\
\Gamma_7 
&\leq Ch^{2k}\dta t \big(h^2+ \|\b^0\|_{k+1}^2\big) \big(h^2 +\|\u^0\|_{k+1}^2\big) \|\b^1\|_1^2 +Ch^{2k}\dta t\|\b^1\|_{k+1}^2 +\frac{1}{2T}\dta t|e_q^1|^2.
\end{aligned} \label{1c_estimateeq4}
\end{align}
For the estimate of the term $\Gamma_8$, by the properties of $L^2$-projection, we obtain
\begin{align}
\begin{aligned}
\Gamma_8 
&= \dta t \big(  (\f^1 - \fh^1,2e_{\u}^1)  + (\g^1 - \gh^1,2e_{\b}^1)\big)\\
&= \dta t \sum_{E\in \Th}\left((\f^1 - \Pi_k^{0,E}\f^1,2e_{\u}^1)_E + (\g^1 - \Pi_k^{0,E}\g^1,2e_{\b}^1)_E \right) \\
&= \dta t \sum_{E\in \Th}\left(((\I-\Pi_k^{0,E})\f^1,(\I - \Pi_0^{0,E})2e_{\u}^1)_E + ((\I-\Pi_k^{0,E})\g^1,(\I - \Pi_0^{0,E})2e_{\b}^1)_E \right) \\
& \leq C  h^k\dta t\big(\|\f^1\|_{k-1}\|\nabla e_{\u}^1\|_0 + \|\g^1\|_{k-1}\|e_{\b}^1\|_1 \big)\\
& \leq Ch^{2k}\dta t \big(\|\f^1\|_{k-1}^2 + \|\g^1\|_{k-1}^2 \big) + \frac{1}{16}\dta t a_{0h}(e_{\u}^1,e_{\u}^1) + \frac{1}{12}\dta t a_{1h}(e_{\b}^1,e_{\b}^1).
\end{aligned} \label{1c_estimateeq5}
\end{align}
According to Lemma \ref{G_estimate}, we can conclude that 
\begin{align}
\Gamma_9 \leq  C (\dta t)^4 +\frac{1}{4}\big(m_{0h}(e_{\u}^1,e_{\u}^1) + m_{1h}(e_{\b}^1,e_{\b}^1) + |e_q^1|^2\big). \label{1c_estimateeq6}
\end{align} 
Substituting \refe{1c_estimateeq2}-\refe{1c_estimateeq6} into \refe{1c_estimateeq1} and using the stability in Lemma \ref{stalility-consistency}, we can obtain the first result of Theorem \ref{1c-error-estimate}.  Then, by triangle inequality and Lemmas \ref{interpolation-DFCVEM}-\ref{interpolation-CVEM}, the second one can be deduced.  
\end{proof}

\begin{theorem} \label{nc-error-estimate}
Let Assumption \ref{regassume} hold. When $h$ and $\dta t$ are small enough, for all $ 1 \leq m \leq N $, the following estimate holds
\begin{align*}
&\|\u^{m+1} - \uh^{m+1}\|_0^2+ \|\b^{m+1} - \bh^{m+1}\|_0^2 + |q^{m+1} - \qh^{m+1}|^2 + \dta t \sum_{n=1}^{m}\Big( \|\nabla(\u^{n+1} - \uh^{n+1})\|_0^2 \\
&+ \|\b^{n+1} - \bh^{n+1}\|_1^2 + \frac{1}{T} |q^{n+1} - \qh^{n+1}|^2\Big)  
\leq C\big((\dta t)^4 + h^{2k} \big).
\end{align*}
\end{theorem}
\begin{proof}
We proceed in several steps for the proof.

Step 1. From \refe{eqMHDeq1} and \refe{fdisMHDeq1}, we get the error equation with the test function $2e_{\u}^{n+1} \in \Zh$
\begin{align}
\begin{aligned}
&m_0(\u_t^{n+1},2e_{\u}^{n+1}) + a_0(\u^{n+1},2e_{\u}^{n+1}) - c_{0}(\u^{n+1},\u^{n+1},2e_{\u}^{n+1}) - c_{1}(\b^{n+1},\b^{n+1},2e_{\u}^{n+1}) \\
&  - m_{0h}(\frac{3\uh^{n+1} - 4\uh^n + \uh^{n-1}}{2\dta t},2e_{\u}^{n+1})  - a_{0h}(\uh^{n+1},2e_{\u}^{n+1})+ \frac{q_h^{n+1}}{q^{n+1}}c_{0h}(\huh^n,\huh^n,2e_{\u}^{n+1})  \\
& + \frac{q_h^{n+1}}{q^{n+1}}c_{1h}(\hbh^n,\hbh^n,2e_{\u}^{n+1}) = (\f^{n+1} - \fh^{n+1},2e_{\u}^{n+1}).
\end{aligned} \label{step1eq1}
\end{align}
From \refe{eqMHDeq2} and \refe{fdisMHDeq2}, taking the test function $2e_{\b}^{n+1}$, we get the error equation
\begin{align}
\begin{aligned}
&m_1(\b_t^{n+1},2e_{\b}^{n+1}) + a_1(\b^{n+1},2e_{\b}^{n+1}) + c_1(2e_{\b}^{n+1},\b^{n+1},\u^{n+1}) \\
&- m_{1h}(\frac{3\b_h^{n+1} - 4\b_h^n + \bh^{n-1}}{2\dta t},2e_{\b}^{n+1}) -a_{1h}(\b_h^{n+1},2e_{\b}^{n+1}) - \frac{q_h^{n+1}}{q^{n+1}}c_{1h}(2e_{\b}^{n+1},\hbh^n,\huh^n) \\
&= (\g^{n+1} - \g_h^{n+1},2e_{\b}^{n+1}).
\end{aligned}  \label{step1eq2}
\end{align}
Adding \refe{step1eq1} and \refe{step1eq2}, we deduce
\begin{align}
\begin{aligned}
&m_{0h}(\frac{3e_{\u}^{n+1} - 4e_{\u}^n + e_{\u}^{n-1}}{2\dta t},2e_{\u}^{n+1}) + m_{1h}(\frac{3e_{\b}^{n+1} - 4e_{\b}^n + e_{\b}^{n-1}}{2\dta t},2e_{\b}^{n+1}) + a_{0h}(e_{\u}^{n+1},2e_{\u}^{n+1}) \\
&+a_{1h}(e_{\b}^{n+1},2e_{\b}^{n+1}) \\
&= \Big\{ m_{0h}(\frac{3\u_I^{n+1} - 4\u_I^n + \u_I^{n-1}}{2\dta t},2e_{\u}^{n+1}) - m_0(\u_t^{n+1},2e_{\u}^{n+1}) + a_{0h}(\u_I^{n+1},2e_{\u}^{n+1}) \\
& \quad- a_0(\u^{n+1},2e_{\u}^{n+1})\Big\}   + \Big\{ m_{1h}(\frac{3\b_I^{n+1} - 4\b_I^n + \b_I^{n-1}}{2\dta t},2e_{\b}^{n+1}) - m_1(\b_t^{n+1},2e_{\b}^{n+1}) \\
& \quad+ a_{1h}(\b_I^{n+1},2e_{\b}^{n+1}) - a_1(\b^{n+1},2e_{\b}^{n+1})\Big\}   +\Big\{c_{0}(\u^{n+1},\u^{n+1},2e_{\u}^{n+1}) \\
& \quad- \frac{q_h^{n+1}}{q^{n+1}}c_{0h}(\huh^n,\huh^n,2e_{\u}^{n+1})\Big\}  + \Big\{c_{1}(\b^{n+1},\b^{n+1},2e_{\u}^{n+1})  - \frac{q_h^{n+1}}{q^{n+1}}c_{1h}(\hbh^n,\hbh^n,2e_{\u}^{n+1})\Big\} \\
& \quad+ \Big\{\frac{q_h^{n+1}}{q^{n+1}}c_{1h}(2e_{\b}^{n+1},\hbh^n,\huh^n) - c_1(2e_{\b}^{n+1},\b^{n+1},\u^{n+1})\Big\}  + \big\{(\f^{n+1} - \fh^{n+1},2e_{\u}^{n+1}) \\
& \quad + (\g^{n+1} - \g_h^{n+1},2e_{\b}^{n+1})\big\}\\
&= \Theta_1 + \Theta_2 + \Theta_3 + \Theta_4 + \Theta_5 + \Theta_6.
\end{aligned}  \label{step1eq3}
\end{align}

Let $\alpha_n = \max_{E \in \Th}\|\Pi_k^{0,E}\huh^n\|_{L^{\infty}(E)}$ and $\beta_n = \max_{E \in \Th} \|\Pi_k^{0,E}\hbh^n\|_{L^{\infty}(E)}$. Using the inverse estimates of polynomial \cite{FEM2008}, the properties of $L^2$-projection, Lemma \ref{stabilityestimate} and stablity in Lemma \ref{stalility-consistency}, then for all $1\leq n \leq m$, there holds
\begin{align}
\begin{aligned}
\dta t \sum_{i=1}^n (\alpha_i^2 + \beta_i^2) 
&= \dta t \sum_{i=1}^n \big( \max_{E \in \Th}\|\Pi_k^{0,E}\huh^i\|_{L^{\infty}(E)}^2 + \max_{E \in \Th} \|\Pi_k^{0,E}\hbh^i\|_{L^{\infty}(E)}^2 \big) \\
& \leq C\dta t \sum_{i=1}^n \big( \max_{E \in \Th} |\Pi_k^{0,E}\huh^i |_{1,E}^2 + \max_{E \in \Th} \|\Pi_k^{0,E}\hbh^i\|_{1,E}^2 \big) \\
& \leq C\dta t \sum_{i=1}^n \big( |\huh^i|_1^2 +  \|\hbh^i\|_1^2 \big) \leq C.
\end{aligned} \label{assume-for-estimate-eq1}
\end{align}
Firstly, making use of the consistency in Lemma \ref{stalility-consistency}, we derive
\begin{align}
\Theta_1
&= m_{0h}(\frac{(3\u_I^{n+1} - 4\u_I^n +\u_I^{n-1}) - (3\u_{\pi}^{n+1} - 4\u_{\pi}^n + \u_{\pi}^{n-1})}{2\dta t},2e_{\u}^{n+1}) \notag\\
&\quad+ \sum_{E\in \Th} m_0^E(\frac{ (3\u_{\pi}^{n+1} - 4\u_{\pi}^n + \u_{\pi}^{n-1}) - (3\u^{n+1} - 4\u^n + \u^{n-1})}{2\dta t},2e_{\u}^{n+1}) \notag\\
& \quad + m_0(\frac{3\u^{n+1} - 4\u^n + \u^{n-1}}{2\dta t} - \u_t^{n+1},2e_{\u}^{n+1} ) + a_{0h}(\u_I^{n+1} - \u_{\pi}^{n+1},2e_{\u}^{n+1}) \notag\\
&\quad + \sum_{E\in \Th} a_0^E(\u_{\pi}^{n+1} - \u^{n+1},2e_{\u}^{n+1}) . \label{step1eq4}
\end{align}
Next, we estimate the right-hand side of \refe{step1eq4} term by term. Using Cauchy-Schwarz inequality, the stability of Lemma \ref{stalility-consistency}, triangle inequality, Lemma \ref{interpolation-DFCVEM} and Lemma \ref{pi-estimate} yields 
\begin{align}
& m_{0h}(\frac{(3\u_I^{n+1} - 4\u_I^n + \u_I^{n-1}) - (3\u_{\pi}^{n+1} - 4\u_{\pi}^n + \u_{\pi}^{n-1})}{2\dta t},2e_{\u}^{n+1}) \notag\\
& \leq \frac{C}{\dta t} \|(3\u_I^{n+1} - 4\u_I^n+ \u_I^{n-1}) - (3\u_{\pi}^{n+1} - 4\u_{\pi}^n + \u_{\pi}^{n-1})\|_0\|e_{\u}^{n+1}\|_0 \notag\\
& \leq \frac{C}{(\dta t)^2}h^{2k+2}|3\u^{n+1} -4\u^n + \u^{n-1}|_{k+1}^2  + \kappa_{0*} \|e_{\u}^{n+1}\|_0^2 \notag\\
& \leq \frac{C}{\dta t}h^{2k+2}\int_{t^{n-1}}^{t^{n+1}} \|\u_{t}\|_{k+1}^2 \,dt +  m_{0h}(e_{\u}^{n+1},e_{\u}^{n+1}). \label{step1eq5}
\end{align}
Applying the similar argument as \refe{step1eq5}, we have
\begin{align}
\begin{aligned}
& \sum_{E\in \Th} m_0^E(\frac{ (3\u_{\pi}^{n+1} - 4\u_{\pi}^n + \u_{\pi}^{n-1}) - (3\u^{n+1} - 4\u^n + \u^{n-1})}{2\dta t},2e_{\u}^{n+1})  \\
& \leq \frac{C}{\dta t}h^{2k+2} \int_{t^{n-1}}^{t^{n+1}}\|\u_t\|_{k+1}^2\,dt + m_{0h}(e_{\u}^{n+1},e_{\u}^{n+1}), \label{step1eq6}
\end{aligned}
\end{align}
and
\begin{align}
\begin{aligned}
 m_0(\frac{3\u^{n+1} - 4\u^n + \u^{n-1}}{2\dta t} - \u_t^{n+1},2e_{\u}^{n+1})  \leq  C(\dta t)^3 \int_{t^{n-1}}^{t^{n+1}}\|\u_{ttt}\|_0^2 \,dt + m_{0h}(e_{\u}^{n+1},e_{\u}^{n+1}). \label{step1eq7}
\end{aligned}
\end{align}
For the estimate of the last two terms, we use Cauchy-Schwarz inequality, the stability of Lemma \ref{stalility-consistency}, triangle inequality, Lemma \ref{interpolation-DFCVEM} and Lemma \ref{pi-estimate} to get
\begin{align}
\begin{aligned}
& a_{0h}(\u_I^{n+1} - \u_{\pi}^{n+1},2e_{\u}^{n+1}) + \sum_{E\in\Th} a_0^E(\u_{\pi}^{n+1} - \u^{n+1},2e_{\u}^{n+1})  \\
& \leq C \sum_{E \in \Th} |\u_I^{n+1} - \u_{\pi}^{n+1}|_{1,E}|e_{\u}^{n+1}|_{1,E} + 2\nu \sum_{E \in \Th} |\u^{n+1} - \u_{\pi}^{n+1}|_{1,E}|e_{\u}^{n+1}|_{1,E} \\
& \leq Ch^{2k} \|\u^{n+1}\|_{k+1}^2  + \frac{1}{6}a_{0h}(e_{\u}^{n+1},e_{\u}^{n+1}). \label{step1eq8}
\end{aligned}
\end{align}
Thus, the estimate of the term $\Theta_1$ can be obtained by combining \refe{step1eq5}-\refe{step1eq8}
\begin{align}
\begin{aligned}
\Theta_1 
&\leq \frac{C}{\dta t}h^{2k+2} \int_{t^{n-1}}^{t^{n+1}}\|\u_t\|_{k+1}^2\,dt + C(\dta t)^3 \int_{t^{n-1}}^{t^{n+1}} \|\u_{ttt}\|_0^2 \,dt  \\
& \quad + Ch^{2k} \|\u^{n+1}\|_{k+1}^2  + 3m_{0h}(e_{\u}^{n+1},e_{\u}^{n+1}) + \frac{1}{6}a_{0h}(e_{\u}^{n+1},e_{\u}^{n+1}). \label{step1eq9}
\end{aligned}
\end{align}
Analogously, we know that
\begin{align}
\begin{aligned}
 \Theta_2  
& \leq \frac{C}{\dta t}h^{2k+2} \int_{t^{n-1}}^{t^{n+1}}\|\b_t\|_{k+1}^2\,dt + C(\dta t)^3 \int_{t^{n-1}}^{t^{n+1}} \|\b_{ttt}\|_0^2 \,dt  \\
& \quad + Ch^{2k} \|\b^{n+1}\|_{k+1}^2 + 3m_{1h}(e_{\b}^{n+1},e_{\b}^{n+1}) + \frac{1}{4}a_{1h}(e_{\b}^{n+1},e_{\b}^{n+1}).  
\end{aligned} \label{step1eq10}
\end{align}
 
Step 2. We observe that  handling the nonlinear terms in \refe{step1eq3} is crucial for the error analysis. In this step, we will provide the estimates of both the nonlinear terms and the load terms. Firstly, we estimate the term $\Theta_3$, which can be split as
\begin{align}
\Theta_3 
& =  \big\{c_0(\u^{n+1},\u^{n+1},2e_{\u}^{n+1}) - c_0(\hu^n,\hu^n,2e_{\u}^{n+1})\big\} + \big\{c_0(\hu^n,\hu^n,2e_{\u}^{n+1}) - c_{0h}(\huh^n,\huh^n,2e_{\u}^{n+1}) \big\} \notag\\
& \quad+ \frac{1}{q^{n+1}}e_q^{n+1}c_{0h}(\huh^n,\huh^n,2e_{\u}^{n+1})\notag\\
& = \Theta_{31} + \Theta_{32} + \Theta_{33} .  \label{step2-1eq1}
\end{align}
We only focus on the terms $\Theta_{31}$ and $\Theta_{32}$ in \refe{step2-1eq1}, the last term will be balanced at the end of the proof. Using the stability of Lemma \ref{stalility-consistency} yields
\begin{align}
\begin{aligned}
\Theta_{31}
&=  c_0(\u^{n+1}-\hu^n,\u^{n+1},2e_{\u}^{n+1}) + c_0(\hu^n,\u^{n+1}-\hu^n,2e_{\u}^{n+1})\\
&\leq  C \big(\|\hu^n\|_1^2 +\|\u^{n+1}\|_1^2 \big)(\dta t)^3\int_{t^{n-1}}^{t^{n+1}} \|\u_{tt}\|_1^2\,dt + \frac{1}{6} a_{0h}(e_{\u}^{n+1},e_{\u}^{n+1}) .  \label{step2-1eq2}
\end{aligned}
\end{align}
Utilizing \refe{nolineareq1} and \refe {nolineareq4}, we conclude that
\begin{align}
\begin{aligned}
\Theta_{32} 
&= c_0(\hu^n,\hu^n,2e_{\u}^{n+1})  - c_{0h}(\hu^n,\hu^n,2e_{\u}^{n+1}) + c_{0h}(\hu^n,\hu^n,2e_{\u}^{n+1})  - c_{0h}(\huh^n,\huh^n,2e_{\u}^{n+1}) \\
& \leq Ch^{2k}\|\hu^n\|_{k+1}^4  +  C\big(1 + \alpha_n^2\big)\big(h^{2k+2}\|\hu^n\|_{k+1}^2 + m_{0h}(\he_{\u}^n,\he_{\u}^n)\big) + \frac{1}{6}a_{0h}(e_{\u}^{n+1},e_{\u}^{n+1}).
\end{aligned} \label{step2-1eq3}
\end{align}
Thus, collecting \refe{step2-1eq2} and \refe{step2-1eq3}, the bound of the term $\Theta_3$ is derived
\begin{align}
\begin{aligned}
\Theta_3
&\leq   C \big(\|\hu^n\|_1^2 +\|\u^{n+1}\|_1^2 \big)(\dta t)^3\int_{t^{n-1}}^{t^{n+1}} \|\u_{tt}\|_1^2\,dt  + \frac{2}{6} a_{0h}(e_{\u}^{n+1},e_{\u}^{n+1}) + \Theta_{33}  \\
&\quad + Ch^{2k}\|\hu^n\|_{k+1}^4  +  C \big(1 + \alpha_n^2 \big) \big(h^{2k+2}\|\hu^n\|_{k+1}^2 + m_{0h}(\he_{\u}^n,\he_{\u}^n)\big) . \label{step2-1eq4}
\end{aligned}
\end{align}
Secondly, in order to estimate the term $\Theta_4$, we rewrite it as 
\begin{align}
\Theta_4 
& = \big\{c_1(\b^{n+1},\b^{n+1},2e_{\u}^{n+1}) - c_1(\hb^n,\hb^n,2e_{\u}^{n+1})\big\}  + \big\{ c_1(\hb^n,\hb^n,2e_{\u}^{n+1}) - c_{1h}(\hbh^n,\hbh^n,2e_{\u}^{n+1})\big\} \notag\\
& \quad +  \frac{1}{q^{n+1}}e_q^{n+1}c_{1h}(\hbh^n,\hbh^n,2e_{\u}^{n+1}) \notag\\
& = \Theta_{41} + \Theta_{42} + \Theta_{43}.  \label{step2-2eq1}
\end{align}
For the estimate of the term $\Theta_{41}$, there holds
\begin{align}
\begin{aligned}
\Theta_{41} 
& =  c_1(\b^{n+1} - \hb^n,\b^{n+1},2e_{\u}^{n+1}) + c_1(\hb^n,\b^{n+1} - \hb^n,2e_{\u} ^{n+1})  \\
& \leq C  \big(\|\hb^n\|_1^2 + \|\b^{n+1}\|_1^2\big)(\dta t)^3 \int_{t^{n-1}}^{t^{n+1}}  \|\b_{tt}\|_1^2\,dt + \frac{1}{6}a_{0h}(e_{\u}^{n+1},e_{\u}^{n+1}).
\end{aligned} \label{step2-2eq2}
\end{align}
From \refe{nolineareq2} and \refe{nolineareq5}, we observe that
\begin{align}
\begin{aligned}
\Theta_{42} 
& =  c_1(\hb^n,\hb^n,2e_{\u}^{n+1}) - c_{1h}(\hb^n,\hb^n,2e_{\u}^{n+1}) +  c_{1h}(\hb^n,\hb^n,2e_{\u}^{n+1}) - c_{1h}(\hbh^n,\hbh^n,2e_{\u}^{n+1}) \\
&= C h^{2k}\big(1 +h^2\|\hb^n\|_2^2 +\|\hb^n\|_k^2\big)\|\hb^n\|_{k+1}^2 + C\|\hb^n\|_2^2m_{1h}(\he_{\b}^n,\he_{\b}^n) +C \beta_n^2 m_{0h}(e_{\u}^{n+1},e_{\u}^{n+1}) \\
& \quad + \frac{1}{40}a_{1h}(\he_{\b}^n,\he_{\b}^n) + \frac{1}{6}a_{0h}(e_{\u}^{n+1},e_{\u}^{n+1}).
\end{aligned} \label{step2-2eq3}
\end{align}
Plugging \refe{step2-2eq2} and \refe{step2-2eq3} into \refe{step2-2eq1}, we obtain
\begin{align}
\Theta_4  
&\leq C \big(\|\hb^n\|_1^2 + \|\b^{n+1}\|_1^2 \big)(\dta t)^3\int_{t^{n-1}}^{t^{n+1}} \|\b_{tt}\|_1^2 \,dt + C h^{2k}\big(1 +h^2\|\hb^n\|_2^2+\|\hb^n\|_k^2 \big)\|\hb^n\|_{k+1}^2 + \Theta_{43} \notag\\
& \quad  + C\|\hb^n\|_2^2m_{1h}(\he_{\b}^n,\he_{\b}^n) + C \beta_n^2 m_{0h}(e_{\u}^{n+1},e_{\u}^{n+1}) + \frac{1}{40}a_{1h}(\he_{\b}^n,\he_{\b}^n) + \frac{2}{6}a_{0h}(e_{\u}^{n+1},e_{\u}^{n+1}) . \label{step2-2eq4}
\end{align}
Thirdly, it is easy to know that the term $\Theta_5$ can be decomposed into the following three terms 
\begin{align}
\Theta_5 
& = \big\{c_1(2e_{\b}^{n+1},\hb^n,\hu^n) - c_1(2e_{\b}^{n+1},\b^{n+1},\u^{n+1})\big\} + \big\{ c_{1h}(2e_{\b}^{n+1},\hbh^n,\huh^n) - c_1(2e_{\b}^{n+1},\hb^n,\hu^n)\big\} \notag\\
& \quad - \frac{1}{q^{n+1}}e_q^{n+1}c_{1h}(2e_{\b}^{n+1},\hbh^n,\huh^n)  \notag\\
&= \Theta_{51} + \Theta_{52} + \Theta_{53}.  \label{step2-3eq1}
\end{align}
By using the stability of Lemma \ref{stalility-consistency}, the following estimate holds 
\begin{align}
\begin{aligned}
\Theta_{51} 
& = c_1(2e_{\b}^{n+1},\hb^n - \b^{n+1},\hu^n) + c_1(2e_{\b}^{n+1},\b^{n+1},\hu^n - \u^{n+1}) \\
& \leq C \big(\|\hu^n\|_1^2 + \|\b^{n+1}\|_1^2  \big)(\dta t)^3\int_{t^{n-1}}^{t^{n+1}} \|\u_{tt}\|_1^2 + \|\b_{tt}\|_1^2 \,dt  + \frac{1}{4}a_{1h}(e_{\b}^{n+1},e_{\b}^{n+1}).
\end{aligned}  \label{step2-3eq2}
\end{align}
From \refe{nolineareq3} and \refe{nolineareq6}, we arrive at 
\begin{align}
\begin{aligned}
\Theta_{52} 
&= c_{1h}(2e_{\b}^{n+1},\hbh^n,\huh^n) - c_{1h}(2e_{\b}^{n+1},\hb^n,\hu^n) +  c_{1h}(2e_{\b}^{n+1},\hb^n,\hu^n) - c_1(2e_{\b}^{n+1},\hb^n,\hu^n) \\
&\leq Ch^{2k}\big(h^2 + \|\hu^n\|_{k+1}^2\big)\|\hb^n\|_{k+1}^2 + C\beta_n^2h^{2k+2}\|\hu^n\|_{k+1}^2 +  C \beta_n^2 m_{0h}(\he_{\u}^n,\he_{\u}^n) \\
& \quad + C m_{1h}(\he_{\b}^n,\he_{\b}^n)  + \frac{1}{4}a_{1h}(e_{\b}^{n+1},e_{\b}^{n+1}).
\end{aligned}  \label{step2-3eq3}
\end{align}
Inserting \refe{step2-3eq2} and \refe{step2-3eq3} into \refe{step2-3eq1}, there holds
\begin{align}
\begin{aligned}
\Theta_5 &\leq C \big(\|\hu^n\|_1^2 + \|\b^{n+1}\|_1^2  \big)(\dta t)^3\int_{t^{n-1}}^{t^{n+1}} \|\u_{tt}\|_1^2 + \|\b_{tt}\|_1^2 \,dt + \frac{2}{4}a_{1h}(e_{\b}^{n+1},e_{\b}^{n+1}) \\
& \quad + Ch^{2k}\big(h^2 + \|\hu^n\|_{k+1}^2\big)\|\hb^n\|_{k+1}^2 + C\beta_n^2h^{2k+2}\|\hu^n\|_{k+1}^2 + C \beta_n^2 m_{0h}(\he_{\u}^n,\he_{\u}^n) \\
& \quad + C m_{1h}(\he_{\b}^n,\he_{\b}^n) + \Theta_{53}.
\end{aligned} \label{step2-3eq4}
\end{align}
Similar to \refe{1c_estimateeq5}, we obtain 
\begin{align}
\begin{aligned}
\Theta_6 \leq  Ch^{2k} \|\f^{n+1}\|_{k-1}^2  + C h^{2k} \|\g^{n+1}\|_{k-1}^2 +  \frac{1}{6}a_{0h}(e_{\u}^{n+1},e_{\u}^{n+1})  + \frac{1}{4} a_{1h}(e_{\b}^{n+1},e_{\b}^{n+1}). 
\end{aligned} \label{step2-4eq1}
\end{align}

Step 3. This step gives the estimate of error on the SAV. Subtracting \refe{fdisMHDeq4} from \refe{eqMHDeq4} at $t=t^{n+1}$ yields
\begin{align}
\begin{aligned}
&\frac{3e_q^{n+1} - 4e_q^{n} + e_{q}^{n-1}}{2\dta t} + \frac{e_q^{n+1}}{T} \\
&= \frac{3q^{n+1} - 4q^n + q^{n-1}}{2\dta t}  - q_t^{n+1}- \frac{1}{q^{n+1}}\Big(c_0(\u^{n+1},\u^{n+1},\u^{n+1})  +c_1(\b^{n+1},\b^{n+1},\u^{n+1})  \\
& \quad- c_1(\b^{n+1},\b^{n+1},\u^{n+1})\Big) + \frac{1}{q^{n+1}} \Big(c_{0h}(\huh^n,\huh^n,\uh^{n+1}) + c_{1h}(\hbh^n,\hbh^n,\uh^{n+1}) \\
& \quad- c_{1h}(\bh^{n+1},\hbh^n,\huh^n) \Big). 
\end{aligned} \label{step3eq1}
\end{align}
Multiplying both sides of \refe{step3eq1} by $2e_q^{n+1}$ leads to
\begin{align}
&\frac{1}{2\dta t}\big( |e_q^{n+1}|^2  + |\he_q^{n+1}|^2 - |e_q^{n}|^2 - |\he_q^n|^2 + |e_q^{n+1} - \he_q^n|^2\big)  + \frac{2}{T}|e_q^{n+1}|^2 \notag\\
&= 2e_q^{n+1}\Big(\frac{3q^{n+1} - 4q^n +q^{n-1}}{2\dta t}  - q_t^{n+1}\Big) - \frac{2}{q^{n+1}}e_q^{n+1} \Big\{ c_0(\u^{n+1} - \hu^n,\u^{n+1},\u^{n+1}) \notag\\
&\quad + c_0(\hu^n,\u^{n+1} - \hu^n,\u^{n+1}) + c_0(\hu^n,\hu^n,\u^{n+1}) - c_{0h}(\hu^n,\hu^n,\u^{n+1}) + c_{0h}(\hu^n,\hu^n,\u^{n+1}) \notag\\
& \quad - c_{0h}(\huh^n,\huh^n,\u^{n+1}) +  c_{0h}(\huh^n,\huh^n,\rho_{\u}^{n+1}) \Big\} - \Theta_{33}  -  \frac{2}{q^{n+1}}e_q^{n+1} \Big\{ c_1(\b^{n+1} -\hb^n,\hb^n,\u^{n+1}) \notag\\
& \quad+ c_1(\hb^n,\hb^n,\u^{n+1}) -c_{1h}(\hb^n,\hb^n,\u^{n+1})  + c_{1h}(\hb^n,\hb^n,\u^{n+1}) - c_{1h}(\hbh^n,\hbh^n,\u^{n+1})  \notag\\
&\quad+ c_{1h}(\hbh^n,\hbh^n,\rho_{\u}^{n+1}) \Big\} - \Theta_{43} + \frac{2}{q^{n+1}}e_q^{n+1} \Big\{  c_1(\b^{n+1},\hb^n,\u^{n+1}-\hu^n)   + c_1(\b^{n+1},\hb^n,\hu^n) \notag\\
& \quad- c_{1h}(\b^{n+1},\hb^n,\hu^n) + c_{1h}(\b^{n+1},\hb^n,\hu^n) - c_{1h}(\b^{n+1},\hbh^n,\huh^n) + c_{1h}(\rho_{\b}^{n+1},\hbh^n,\huh^n) \Big\} - \Theta_{53} \notag\\
&= \Theta_{71} + \Theta_{72} - \Theta_{33} + \Theta_{73} - \Theta_{43} + \Theta_{74} - \Theta_{53}. \label{step3eq2}
\end{align} 
In fact, due to the term $\frac{2e_q^{n+1}}{q^{n+1}} c_1(\b^{n+1},\b^{n+1}-\hb^n,\u^{n+1})$ being balanced, the equality \refe{step3eq2} holds. We first give the estimate of the term $\Theta_{71}$
\begin{align}
\begin{aligned}
\Theta_{71}  
 \leq \frac{1}{4T} |e_q^{n+1}|^2 + C(\dta t)^3 \int_{t^{n-1}}^{t^{n+1}}|q_{ttt}|^2 \,dt. 
\end{aligned} \label{step3eq3}
\end{align}
Similar to $\Theta_3,\Theta_4$ and $\Theta_5$, with a few modifications, we have
\begin{align}
\begin{aligned}
\Theta_{72} 
&\leq  C\big(\|\hu^n\|_1^2 + \|\u^{n+1}\|_1^2 \big)\|\u^{n+1}\|_1^2(\dta t)^3\int_{t^{n-1}}^{n+1} \|\u_{tt}\|_1^2\,dt  + Ch^{2k}\|\hu^n\|_{k+1}^4\|\u^{n+1}\|_1^2 \\
& \quad  + C\big(1 + \alpha_n^2 \big)h^{2k+2}\|\hu^{n}\|_{k+1}^2\|\u^{n+1}\|_1^2 + C(1 + \alpha_n^2)\|\u^{n+1}\|_1^2m_{0h}(\he_{\u}^n,\he_{\u}^n) \\
& \quad+ \frac{1}{4T}|e_q^{n+1}|^2 + C\alpha_n^2 h^{2k}\|\u^{n+1}\|_{k+1}^2, \\
\Theta_{73} 
&\leq  C \|\hb^n\|_1^2\|\u^{n+1}\|_1^2(\dta t)^3\int_{t^{n-1}}^{t^{n+1}}\|\b_{tt}\|_1^2\,dt + Ch^{2k}\big(h^2\|\hb^n\|_2^2 + \|\hb^n\|_k^2 \big)\|\hb^n\|_{k+1}^2 \|\u^{n+1}\|_1^2 \\
& \quad + Ch^{2k}\|\hb^n\|_{k+1}^2   + C\|\hb^n\|_2^2\|\u^{n+1}\|_1^2m_{1h}(\he_{\b}^n,\he_{\b}^n) +\frac{1}{40}a_{1h}(\he_{\b}^n,\he_{\b}^n) \\
& \quad  + C\beta_n^2 \big(1 +\|\u^{n+1}\|_0^2\big)|e_q^{n+1}|^2+  \frac{1}{4T}|e_q^{n+1}|^2 + C\|\hbh^n\|_1^2 h^{2k+2}\|\u^{n+1}\|_{k+1}^2,\\
\Theta_{74}
& \leq  C\|\hb^n\|_1^2\|\b^{n+1}\|_1(\dta t)^3\int_{t^{n-1}}^{n+1} \|\u_{tt}\|_1^2\,dt + Ch^{2k}\big(h^2 + \|\hu^n\|_{k+1}^2 \big)\|\hb^n\|_{k+1}^2 \|\b^{n+1}\|_1^2 \\
& \quad + C\beta_n^2 h^{2k+2}\|\hu^n\|_{k+1}^2\|\b^{n+1}\|_1^2 + C\beta_n^2\|\b^{n+1}\|_1^2m_{0h}(\he_{\u}^n,\he_{\u}^n) + C\|\b^{n+1}\|_1^2m_{1h}(\he_{\b}^n,\he_{\b}^n) \\
&\quad + \frac{1}{4T}|e_q^{n+1}|^2  +  C\beta_n^2 h^{2k}\|\b^{n+1}\|_{k+1}^2.
\end{aligned} \label{step3eq4}
\end{align}
Next, plugging \refe{step3eq3} and \refe{step3eq4} into \refe{step3eq2} and using Assumption \ref{regassume}, we deduce 
\begin{align}
&\frac{1}{2\dta t}\Big( |e_q^{n+1}|^2  + |\he_q^{n+1}|^2 - |e_q^{n}|^2 - |\he^n|^2 \Big)  + \frac{1}{T}|e_q^{n+1}|^2 \notag\\
& \leq  C(\dta t)^3\int_{t^{n-1}}^{n+1} \|\u_{tt}\|_1^2 + \|\b_{tt}\|_1^2 + |q_{ttt}|^2\,dt + C(1 + \alpha_n^2 + \beta_n^2)h^{2k}  \notag\\
& \quad  + C\big(1 + \alpha_n^2  + \beta_n^2\big)h^{2k+2} + C\big(1 + \alpha_n^2 + \beta_n^2\big)m_{0h}(\he_{\u}^n,\he_{\u}^n) + C\beta_n^2 |e_q^{n+1}|^2  \notag\\
& \quad + C m_{1h}(\he_{\b}^n,\he_{\b}^n) +   C\|\hbh^n\|_1^2 h^{2k+2} + \frac{1}{40}a_{1h}(\he_{\b}^n,\he_{\b}^n)  - \Theta_{33} - \Theta_{43} - \Theta_{53}. \label{step3eq5}
\end{align}

Step 4. Inserting \refe{step1eq9}, \refe{step1eq10}, \refe{step2-1eq4}, \refe{step2-2eq4}, \refe{step2-3eq4}, \refe{step2-4eq1} into \refe{step1eq3}, then adding \refe{step3eq5} and applying the Assumption \ref{regassume}, we obtain
\begin{align}
&\frac{1}{2}\Big( m_{0h}(e_{\u}^{n+1},e_{\u}^{n+1}) + m_{0h}(\he_{\u}^{n+1},\he_{\u}^{n+1}) + m_{1h}(e_{\b}^{n+1},e_{\b}^{n+1}) + m_{1h}(\he_{\b}^{n+1},\he_{\b}^{n+1}) + |e_q^{n+1}|^2 \notag\\
&+ |\he_{q}^{n+1}|^2 \Big) + \dta t \Big( a_{0h}(e_{\u}^{n+1},e_{\u}^{n+1})  + a_{1h}(e_{\b}^{n+1},e_{\b}^{n+1})  + \frac{1}{T}|e_q^{n+1}|^2\Big)  \notag\\
& \leq  \frac{1}{2}\Big(m_{0h}(e_{\u}^n,e_{\u}^n) +m_{0h}(\he_{\u}^n,\he_{\u}^n)  + m_{1h}(e_{\b}^n,e_{\b}^n) + m_{1h}(\he_{\b}^{n},\he_{\b}^n) 
+ |e_q^n|^2 +|\he_q^n|^2 \Big) \notag\\
& \quad  + C \dta t\big(1 +\beta_n^2\big) m_{0h}(e_{\u}^{n+1},e_{\u}^{n+1})  + C \dta t m_{1h}(e_{\b}^{n+1},e_{\b}^{n+1}) + C \dta t m_{1h}(\he_{\b}^n,\he_{\b}^n) \notag \\
& \quad  + C\dta t\big(1 + \alpha_n^2 + \beta_n^2\big)m_{0h}(\he_{\u}^n,\he_{\u}^n) + C\dta t\beta_n^2 |e_q^{n+1}|^2 + \frac{\dta t}{20} a_{1h}(\he_{\b}^n,\he_{\b}^n)  \notag\\
& \quad + C(\dta t)^4\int_{t^{n-1}}^{t^{n+1}} \|\u_{ttt}\|_0^2 + \|\b_{ttt}\|_0^2 + |q_{ttt}|^2 \,dt  + C(\dta t)^4\int_{t^{n-1}}^{t^{n+1}} \|\u_{tt}\|_1^2 + \|\b_{tt}\|_1^2 \,dt \notag\\
& \quad   + Ch^{2k+2}\int_{t^{n-1}}^{t^{n+1}} \|\u_t\|_{k+1}^2 + \|\b_t\|_{k+1}^2 \,dt + C\dta t(1 + \alpha_n^2 + \beta_n^2)h^{2k}+ C\dta t\|\hbh^n\|_1^2 h^{2k+2} \notag\\
& \quad    + C\dta t\big(1 + \alpha_n^2  + \beta_n^2\big)h^{2k+2}  . \label{step4eq1}
\end{align}
Summing up \refe{step4eq1} from $n=1$ to $m$, utilizing Cauchy-Schwarz inequality,  Lemma \ref{nc-bound}, Assumption \ref{regassume} and \refe{assume-for-estimate-eq1}, we deduce
\begin{align}
&\frac{1}{2}\Big( m_{0h}(e_{\u}^{m+1},e_{\u}^{m+1}) + m_{0h}(\he_{\u}^{m+1},\he_{\u}^{n+1}) + m_{1h}(e_{\b}^{m+1},e_{\b}^{m+1}) + m_{1h}(\he_{\b}^{m+1},\he_{\b}^{m+1}) + |e_q^{m+1}|^2\notag\\
&+ |\he_{q}^{m+1}|^2 \Big) + \dta t\sum_{n=1}^m \Big( \frac{1}{2} a_{0h}(e_{\u}^{n+1},e_{\u}^{n+1})  + a_{1h}(e_{\b}^{n+1},e_{\b}^{n+1})  + \frac{1}{T}|e_q^{n+1}|^2\Big)  \notag\\ 
&\leq \frac{1}{2}\Big(m_{0h}(e_{\u}^1,e_{\u}^1) +m_{0h}(\he_{\u}^1,\he_{\u}^1)  + m_{1h}(e_{\b}^1,e_{\b}^1) + m_{1h}(\he_{\b}^1,\he_{\b}^1) 
+ |e_q^1|^2 +|\he_q^1|^2 \Big) \notag\\
& \quad + \frac{1}{10}\dta t \big(5a_{1h}(e_{\b}^1,e_{\b}^1) + a_{1h}(e_{\b}^0,e_{\b}^0) \big) + C_0\dta t \sum_{n=1}^m \big(1 + \alpha_n^2 + \beta_n^2\big) \Big(  m_{0h}(e_{\u}^{n+1},e_{\u}^{n+1}) \notag\\
& \quad + m_{0h}(\he_{\u}^n,\he_{\u}^n) + m_{1h}(e_{\b}^{n+1},e_{\b}^{n+1}) +  m_{1h}(\he_{\b}^n,\he_{\b}^n) + |e_q^{n+1}|^2 \Big)   + C(\dta t)^4 + Ch^{2k} \big( 1 + h^2\big) ,
\end{align}
where $C_0$ is a positive constant independent of $h$ and $\dta t$. We assume that when $h$ and $\dta t$ are small enough, there is 
\begin{align}
\dta t C_0(1 + \alpha_n^2 + \beta_n^2 ) \leq \frac{1}{2}, \quad \forall\, 1 \leq n \leq m.  \label{assume-for-estimate-eq2}
\end{align}
Then, using the stability of Lemma \ref{stalility-consistency}, Theorem \ref{1c-error-estimate}, Lemma \ref{gronwelleq} and \refe{assume-for-estimate-eq1}, we infer
\begin{align*}
&\frac{1}{2}\big(\|e_{\u}^{m+1}\|_0^2+ \|e_{\b}^{m+1}\|_0^2 + |e_q^{m+1}|^2 \big) + \dta t \sum_{n=1}^{m}\big( \frac{1}{2}\|\nabla e_{\u}^{n+1}\|_0^2 + \|e_{\b}^{n+1}\|_1^2 + \frac{1}{T} |e_q^{n+1}|^2\big)   \\
&\leq C\exp(2C_0 C)\big((\dta t)^4 + h^{2k} \big).
\end{align*}
By triangle inequality and Lemmas \ref{interpolation-DFCVEM}-\ref{interpolation-CVEM}, we drive the result of Theorem \ref{nc-error-estimate}.

Step 5. In this step, we will prove that the assumption \refe{assume-for-estimate-eq2} holds by mathematical induction. First, we show that \refe{assume-for-estimate-eq2} holds for $n=1$. Similar to \refe{assume-for-estimate-eq1}, by using Theorem \ref{1c-error-estimate}, we have
\begin{align*}
\dta t C_0(1 + \alpha_1^2 +  \beta_1^2) 
& = \dta t C_0\big(1 + \max_{E \in \Th}\|\Pi_k^{0,E}\huh^1\|_{L^{\infty}(E)}^2 + \max_{E \in \Th} \|\Pi_k^{0,E}\hbh^1\|_{L^{\infty}(E)}^2\big) \\
& \leq  C_0C\dta t \big(1 + |\huh^1 |_1^2 + \|\hbh^1\|_1^2 \big) \\
& \leq C_0C\dta t \big(1 + |\huh^1 - \hu^1|_1^2 + |\hu^1|_1^2 + \|\hbh^1 - \hb^1\|_1^2 + \|\hb^1\|_1^2 \big) \\
& \leq C_0C \big( \dta t +  (\dta t)^4 +  h^{2k}\big) \leq \frac{1}{2}.
\end{align*}
Suppose \refe{assume-for-estimate-eq2} holds for $n = m$. Then, similar to $n=1$, by the result of Theorem \ref{nc-error-estimate}, we  can prove that \refe{assume-for-estimate-eq2} holds for $n=m+1$. Thus the proof is completed.
\end{proof}

\section{Numerical experiments}
\begin{figure}[H]
\vspace{-0.5cm}
\centering
\subfigure
{
    \begin{minipage}[b]{.3\linewidth}
        \centering
        \includegraphics[scale=0.4]{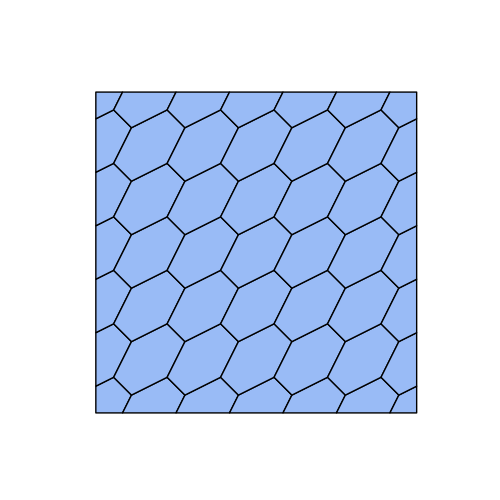}
    \end{minipage}
}
\subfigure
{
 	\begin{minipage}[b]{.3\linewidth}
        \centering
        \includegraphics[scale=0.4]{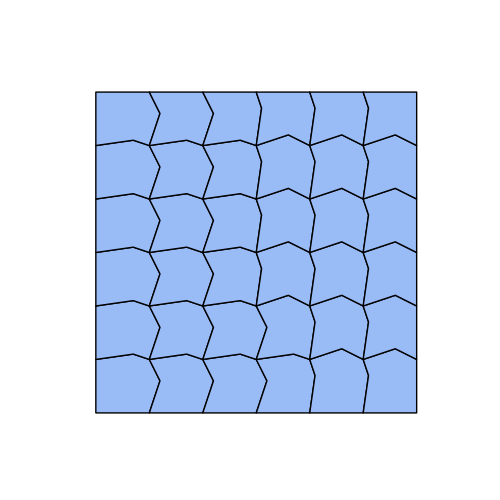}
    \end{minipage}
}
\subfigure
{
 	\begin{minipage}[b]{.3\linewidth}
        \centering
        \includegraphics[scale=0.4]{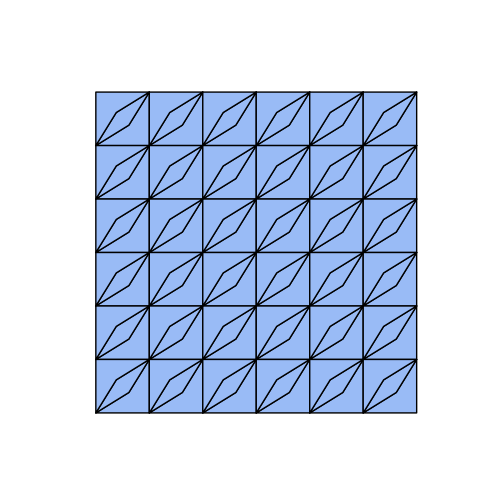}
    \end{minipage}
}
\vspace{-0.5cm}
\caption{Polygonal meshes: $\Th^1$(left), $\Th^2$(middle) and $\Th^3$(right).} \label{Fig1}
\end{figure}

In this section, we present two numerical experiments to validate the theoretical results. The numerical experiments are implemented by Fealpy package \cite{Wei-Fealpy}. Example \ref{example1} and Example \ref{example2} are performed with $k=1, \dta t \approx h$ (i.e. $N+1 = \lceil T/h\rceil$). Example \ref{example1} is used to verify the error estimates of the fully discrete BDF2 IMEX SAV scheme on three types of polygonal meshes (see Figure \ref{Fig1}). Example \ref{example2} is used to simulate a driven cavity flow. 
In order to compute the VEM errors, we consider the following computable error quantities:
\begin{align*}
\error(\u,L^2) &=  \sqrt{\sum_{E\in\Th} \|\u - \Pi_k^{0,E}\uh\|_{0,E}^2}\, ,\qquad \error(\u,H^1) = \sqrt{\sum_{E\in\Th} |\u - \Pi_k^{\nabla,E}\uh|_{1,E}^2}\, ,\\
\error(\b,L^2) &=  \sqrt{\sum_{E\in\Th} \|\b - \Pi_k^{0,E}\bh\|_{0,E}^2}\,,\qquad \error(\b,H^1) = \sqrt{\sum_{E\in\Th} \|\b - \Pi_k^{\nabla,E}\bh\|_{1,E}^2}\,,\\
\error(p,L^2) &= \|p - p_h\|_0, \qquad \error(q) = |q - \qh| .
\end{align*}
\begin{example} \label{example1}
Let the domain $\O = (0,1)^2$, $T = 1$, the equations parameters $\nu=1,\mu=1,\sigma = 1$ and the exact solution be 
\begin{align*}
\u &= \begin{pmatrix}
\pi \exp(-t) \cos(t) \sin(\pi x_1)^2 \sin(\pi x_2)\cos(\pi x_2)\\
-\pi \exp(-t) \cos(t) \sin(\pi x_1) \sin(\pi x_2)^2\cos(\pi x_1)
\end{pmatrix},\\
\b &= \begin{pmatrix}
\exp(-t)\cos(t)\sin(\pi x_1)\cos(\pi x_2)\\
-\exp(-t)\cos(t)\cos(\pi x_1)\sin(\pi x_2)
\end{pmatrix},\\
p &= \exp(-t)\cos(t)\cos(\pi x_1)\cos(\pi x_2).
\end{align*}
\end{example}
In numerical simulation of this example, the above mentioned polygonal meshes are used with $5$ refinement levels (level $i$, $i=1,...,5$). The mesh sizes of different types meshes may be different at the same refinement level. In Tables \ref{table1} and \ref{table2}, we show the specific mesh sizes of refinement levels of each type mesh.  However, for simplicity, we show the refinement levels of meshes in Table \ref{table3}.
\begin{table}[H] 
\centering
\caption{Convergence rates of velocity and pressure at $t=1, \dta t \approx  h$.}
\begin{tabular}{ccclclcl}
\hline
mesh  & $h$  & $\error(\u,L^2)$ & Rate & $\error(\u,H^1)$ & Rate   & $\error(p,L^2)$ & Rate \\ \hline 
$\Th^1$
  & 0.2981 & 8.8977e-02 & -    & 6.1953e-01 & -    & 1.4799e-01 & -    \\      
  & 0.1491 & 1.9938e-02 & 2.16 & 3.1631e-01 & 0.97 & 4.1332e-02 & 1.84 \\  
  & 0.0745 & 5.0060e-03 & 1.99 & 1.6036e-01 & 0.98 & 1.2243e-02 & 1.76 \\    
  & 0.0373 & 1.2585e-03 & 1.99 & 8.0699e-02 & 0.99 & 4.2659e-03 & 1.52 \\     
  & 0.0186 & 3.1585e-04 & 1.99 & 4.0476e-02 & 1.0  & 1.8117e-03 & 1.24 \\ \hline 
$\Th^2$
  & 0.2828 & 7.5109e-02 & -    & 6.3204e-01 & -    & 6.6826e-02 & -    \\    
  & 0.1414 & 1.5820e-02 & 2.25 & 3.1291e-01 & 1.01 & 1.8550e-02 & 1.85 \\   
  & 0.0707 & 3.8417e-03 & 2.04 & 1.5617e-01 & 1.0  & 7.3819e-03 & 1.33 \\   
  & 0.0354 & 9.5227e-04 & 2.01 & 7.8031e-02 & 1.0  & 3.3778e-03 & 1.13 \\ 
  & 0.0177 & 2.3750e-04 & 2.0  & 3.9005e-02 & 1.0  & 1.6410e-03 & 1.04 \\ \hline
$\Th^3$ 
  & 0.2828 & 2.8077e-02 & -    & 5.0891e-01 & -    & 1.5975e-01 & -    \\  
  & 0.1414 & 6.7140e-03 & 2.06 & 2.5939e-01 & 0.97 & 6.0931e-02 & 1.39 \\  
  & 0.0707 & 1.6754e-03 & 2.0  & 1.3019e-01 & 0.99 & 2.5078e-02 & 1.28 \\ 
  & 0.0354 & 4.1805e-04 & 2.0  & 6.5148e-02 & 1.0  & 1.1546e-02 & 1.12 \\        
  & 0.0177 & 1.0440e-04 & 2.0  & 3.2577e-02 & 1.0  & 5.6011e-03 & 1.04 \\ \hline
\end{tabular}  \label{table1}
\end{table} 

\begin{table}[H] 
\centering
\caption{Convergence rates of magnetic field and auxiliary variable at $t=1, \dta t \approx  h$.}
\begin{tabular}{ccclclcl}
\hline
mesh  & $h$  & $\error(\b,L^2)$ & Rate & $\error(\b,H^1)$ & Rate   & $\error(q)$ & Rate \\ \hline 
$\Th^1$
  & 0.2981 & 6.8322e-03 & -    & 1.5115e-01 & -    & 3.2151e-02 & -    \\     
  & 0.1491 & 2.0580e-03 & 1.73 & 8.0318e-02 & 0.91 & 1.2771e-02 & 1.33 \\      
  & 0.0745 & 5.4881e-04 & 1.91 & 4.1238e-02 & 0.96 & 2.9835e-03 & 2.1  \\    
  & 0.0373 & 1.4084e-04 & 1.96 & 2.0869e-02 & 0.98 & 6.5250e-04 & 2.19 \\      
  & 0.0186 & 3.5598e-05 & 1.98 & 1.0494e-02 & 0.99 & 1.4922e-04 & 2.13 \\ \hline 
$\Th^2$
  & 0.2828 & 5.9764e-03 & -    & 1.6444e-01 & -    & 2.3447e-02 & -    \\    
  & 0.1414 & 1.4814e-03 & 2.01 & 8.1240e-02 & 1.02 & 8.5508e-03 & 1.46 \\      
  & 0.0707 & 3.7443e-04 & 1.98 & 4.0566e-02 & 1.0  & 1.8849e-03 & 2.18 \\     
  & 0.0354 & 9.4487e-05 & 1.99 & 2.0275e-02 & 1.0  & 3.9766e-04 & 2.24 \\     
  & 0.0177 & 2.3753e-05 & 1.99 & 1.0136e-02 & 1.0  & 8.7973e-05 & 2.18 \\ \hline
$\Th^3$ 
  & 0.2828 & 9.4651e-03 & -    & 1.3821e-01 & -    & 7.5529e-03 & -    \\     
  & 0.1414 & 2.4971e-03 & 1.92 & 6.9984e-02 & 0.98 & 1.9136e-03 & 1.98 \\       
  & 0.0707 & 6.2964e-04 & 1.99 & 3.5023e-02 & 1.0  & 3.7825e-04 & 2.34 \\     
  & 0.0354 & 1.5783e-04 & 2.0  & 1.7515e-02 & 1.0  & 7.8206e-05 & 2.27 \\     
  & 0.0177 & 3.9493e-05 & 2.0  & 8.7580e-03 & 1.0  & 1.7462e-05 & 2.16 \\ \hline
\end{tabular}\label{table2}
\end{table} 
In Table \ref{table1}, we present errors and convergence rates of the velocity and pressure at final time $T=1$. In Table \ref{table2}, errors and convergence rates of the magnetic field and the auxiliary variable are showed at $T=1$. Table \ref{table3} displays the $L^2$-norm of $\div \uh$ at $T=1$. From Tables \ref{table1} and \ref{table2}, the convergence rates of the velocity, magnetic field and the auxiliary variable are consistency with the results in Theorem \ref{nc-error-estimate}. And the convergence rates of the pressure are satisfying. From Table \ref{table3}, we know that the discrete velocity is divergence-free up to machine precision. 
\begin{table}[H]  
\centering
\caption{The $L^2$-norm of $\div \uh$ at $t=1, \dta t \approx  h$.}
\begin{tabular}{cccccc}
\hline 
mesh    & level 1    & level 2    & level 3    & level 4    & level 5    \\ \hline
$\Th^1$ & 6.1184e-16 & 1.5603e-15 & 4.9325e-15 & 7.8330e-15 & 4.7330e-14 \\  \hline
$\Th^2$ & 9.8100e-16 & 1.8607e-15 & 9.9034e-15 & 7.0958e-15 & 5.1457e-14 \\ \hline
$\Th^3$ & 1.1361e-15 & 2.8907e-15 & 8.1725e-15 & 1.6305e-14 & 7.6550e-14 \\ \hline
\end{tabular} \label{table3}
\end{table} 

\begin{figure}[H]
\vspace{-0.6cm}
\centering
\subfigure
{
    \begin{minipage}[b]{.45\linewidth}
        \centering
        \includegraphics[scale=0.15]{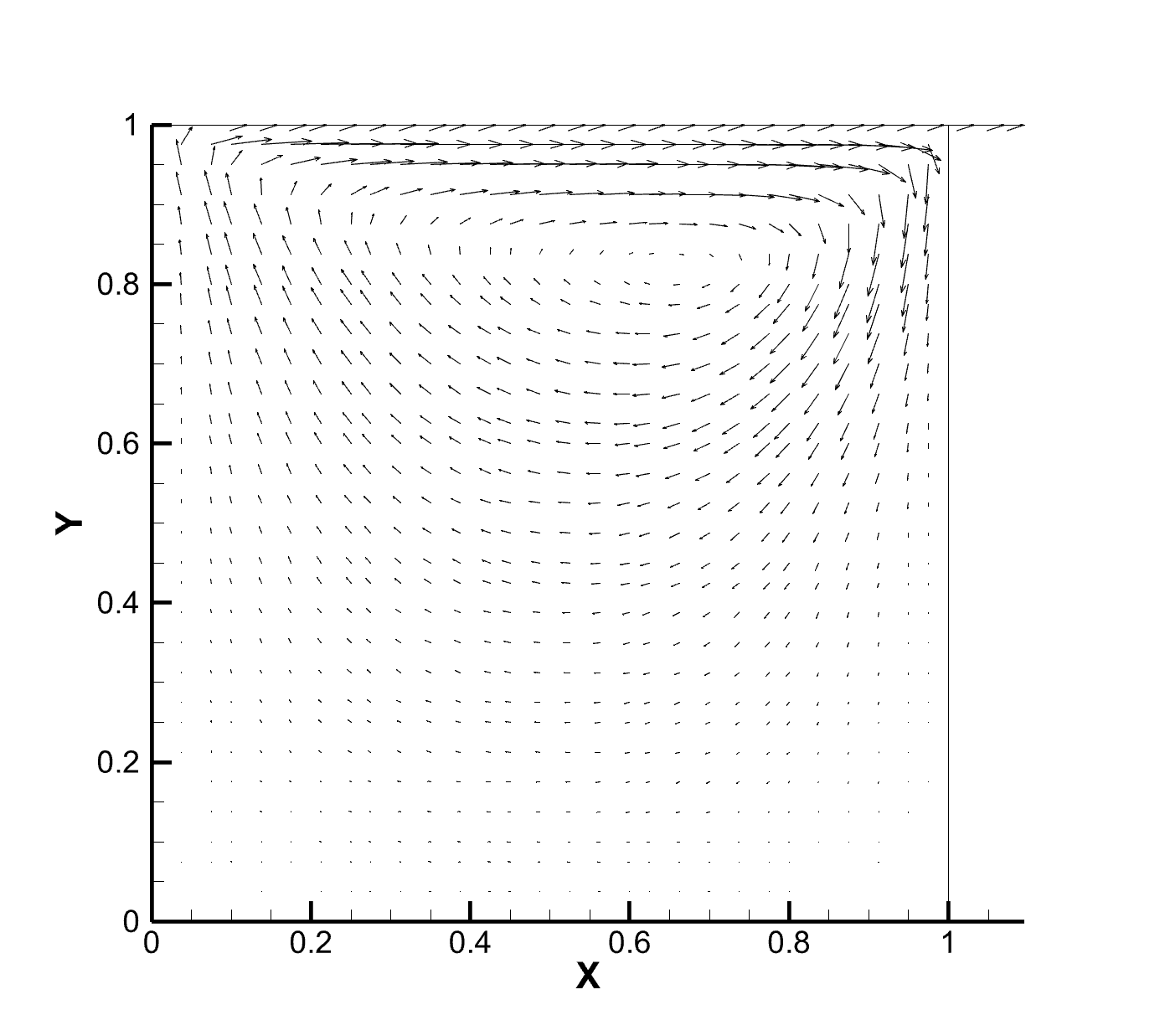}
    \end{minipage}
}
\subfigure
{
 	\begin{minipage}[b]{.45\linewidth}
        \centering
        \includegraphics[scale=0.15]{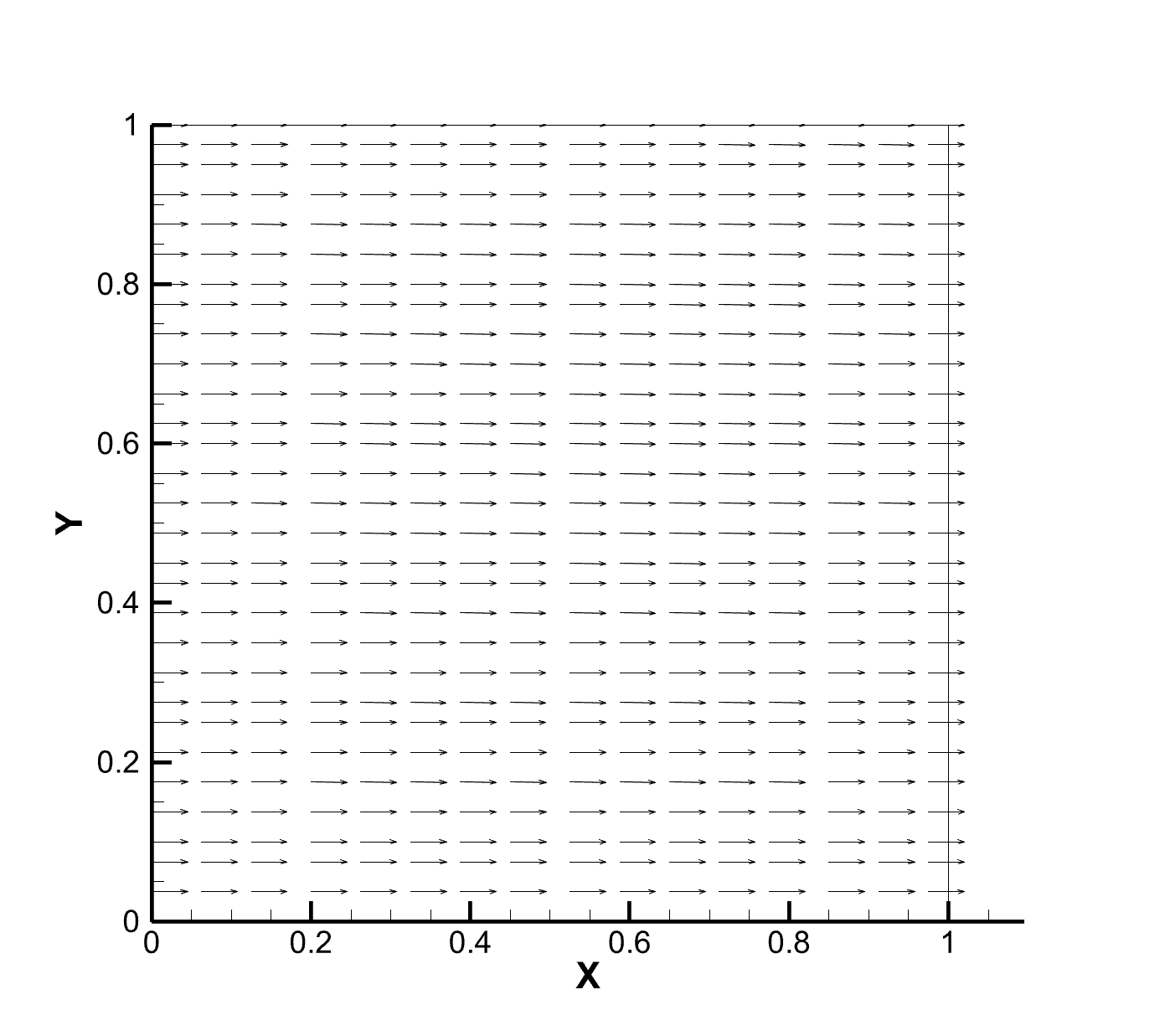}
    \end{minipage}
}
\vspace{-0.5cm}
\caption{The vector graphs of velocity (left) and magnetic (right) field at $t=1$.} \label{Fig2}
\end{figure}

\begin{figure}[H]
\vspace{-0.6cm}
\centering
\subfigure
{
    \begin{minipage}[b]{.45\linewidth}
        \centering
        \includegraphics[scale=0.15]{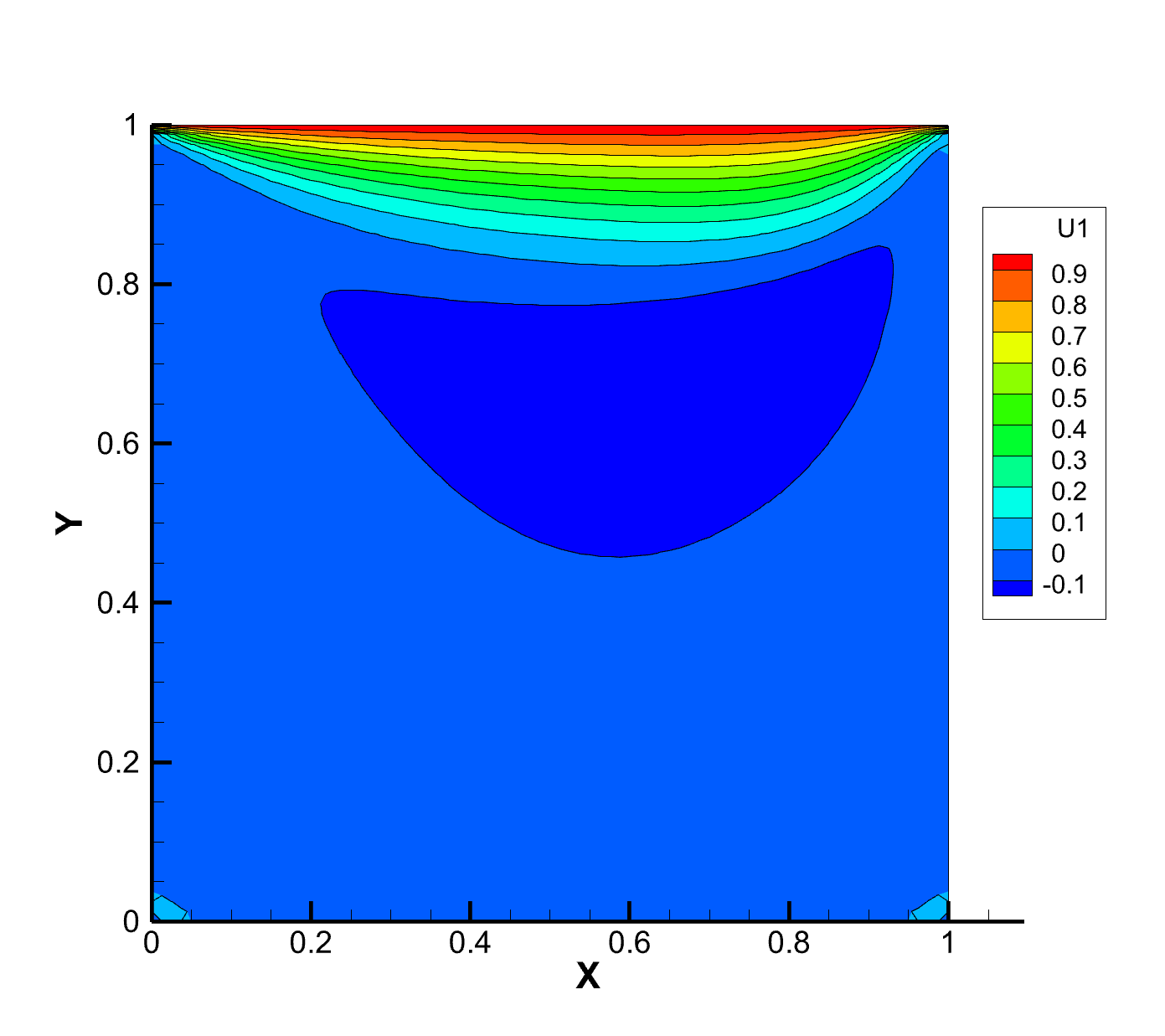}
    \end{minipage}
}
\subfigure
{
 	\begin{minipage}[b]{.45\linewidth}
        \centering
        \includegraphics[scale=0.15]{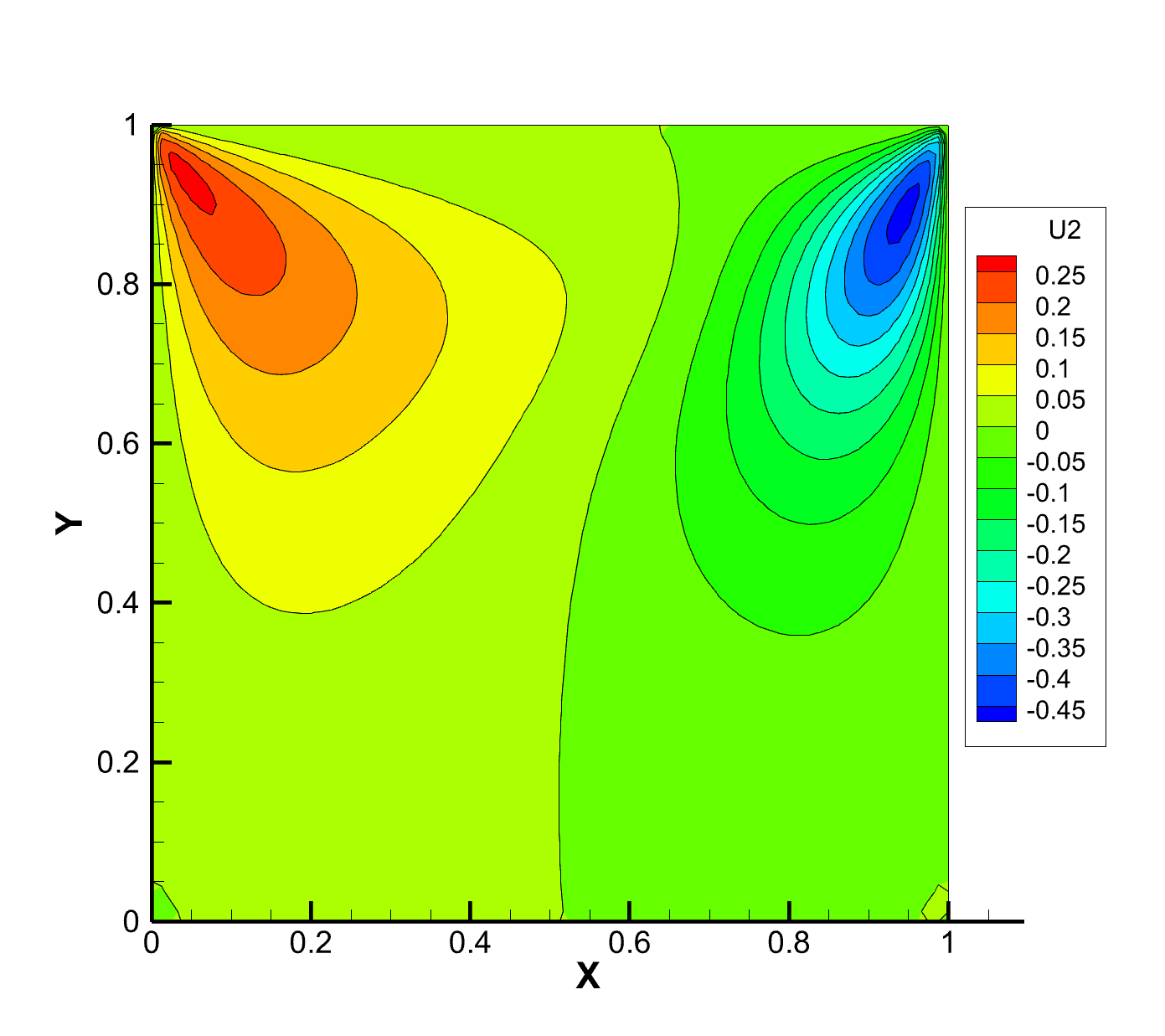}
    \end{minipage}
}
\subfigure
{
    \begin{minipage}[b]{.45\linewidth}
        \centering
        \includegraphics[scale=0.15]{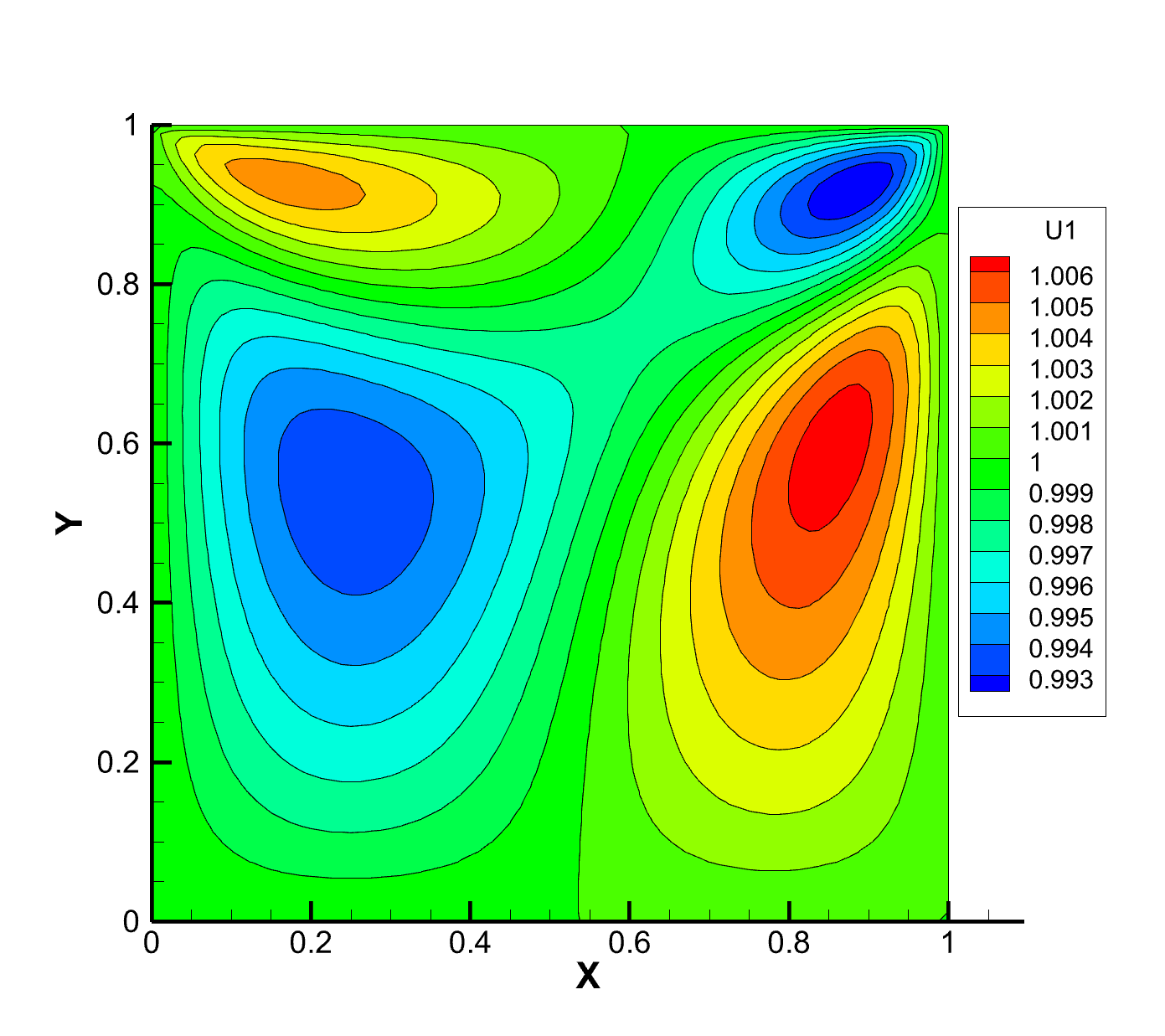}
    \end{minipage}
}
\subfigure
{
 	\begin{minipage}[b]{.45\linewidth}
        \centering
        \includegraphics[scale=0.15]{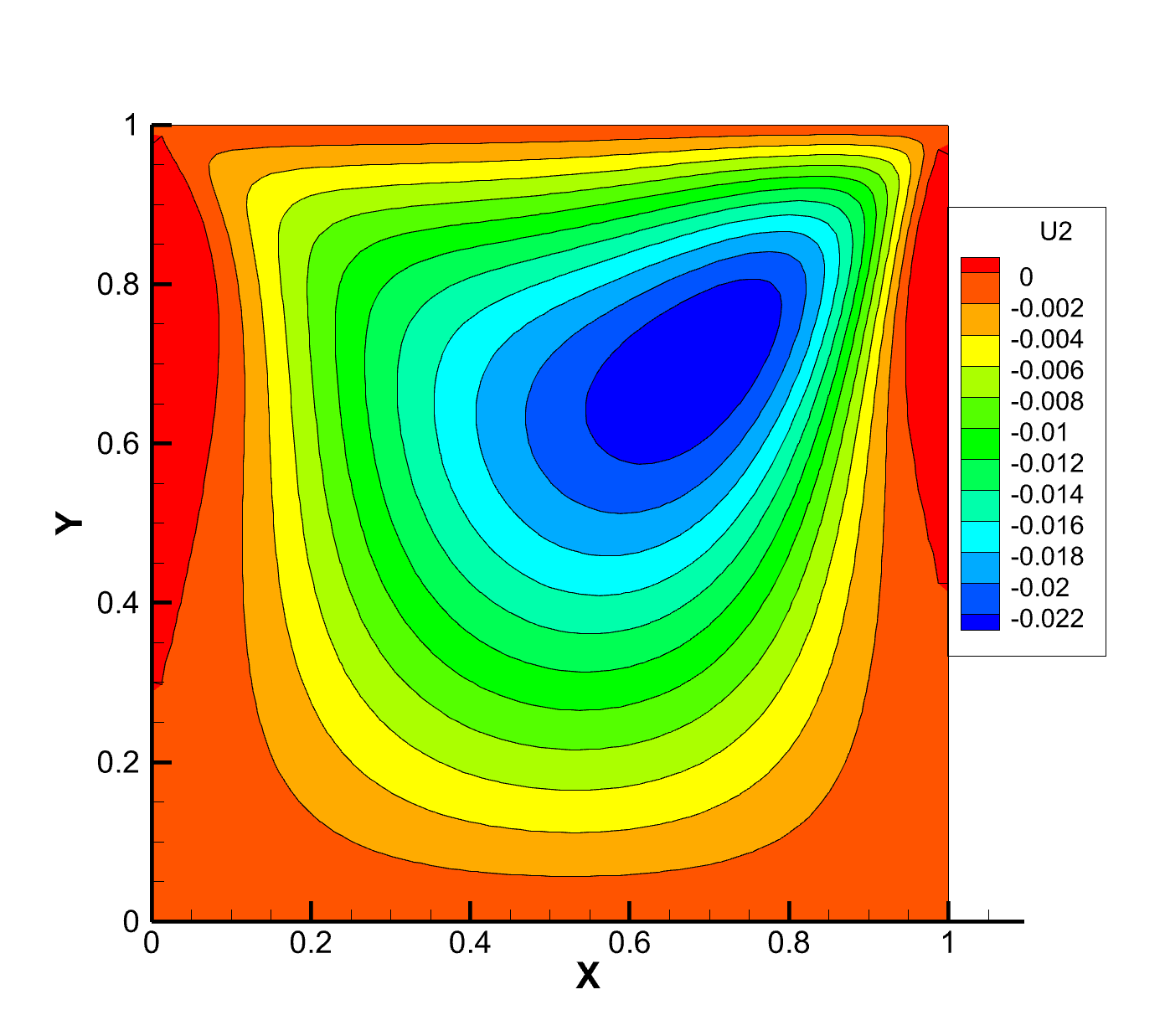}
    \end{minipage}
}
\vspace{-0.5cm}
\caption{The contours of velocity (upper) and magnetic (lower) field at $t=1$.} \label{Fig3}
\end{figure}

\begin{example} \label{example2}
We consider a driven cavity flow, which is a classic benchmark from computational fluid mechanics. Let the domain $\O = (0,1)^2,T=1,\nu = 0.01,\sigma=100,\mu = 1$. The considered boundary conditions are  
\begin{align*}
&\u = \0 \qquad             \mbox{on} \;x_1= \pm 1 \; \mbox{and} \; x_2= 0,\\
&\u = (1,0)^{\mathrm{T}}     \qquad \mbox{on}\; x_2 = 1, \\
&\n \times \b = \n \times (1,0)^{\mathrm{T}}  \qquad \mbox{on} \,\pt \O.
\end{align*}
\end{example}
We perform this example on rectangular meshes with $h =0.0177$. In Figure \ref{Fig2}, we display the vector graphs of the velocity and magnetic field at final time. In Figure \ref{Fig3}, we show the contours of all components of the velocity and magnetic field at final time.

\section*{Acknowledgements}
The research was supported by the National Natural Science Foundation of China (Nos: 12071404, 11971410, 12026254, 12071402, 12261131501), Young Elite Scientist Sponsorship Program by CAST (No: 2020QNRC001), Science and Technology Innovation Program of Hunan Province (No: 2024RC3158), Key Project of Scientific Research Project of Hunan Provincial Department of Education (No: 22A0136), Postgraduate Scientific Research Innovation Project of Hunan Province (No: CX20230614), Project of Scientiﬁc Research Fund of the Hunan Provincial Science and Technology Department (No: 2023GK2029, 2024ZL5017), and Program for Science and Technology Innovative Research Team in Higher Educational Institutions of Hunan Province of China.

\bibliographystyle{abbrv}
\bibliography{references}

\begin{thebibliography}{10}

\bibitem{stream-VEM-NS}
D.~Adak, D.~Mora, S.~Natarajan, and A.~Silgado.
\newblock A virtual element discretization for the time dependent
  {N}avier--{S}tokes equations in stream-function formulation.
\newblock {\em ESAIM: Mathematical Modelling and Numerical Analysis},
  55(5):2535--2566, 2021.

\bibitem{VEM-equivalent}
B.~Ahmad, A.~Alsaedi, F.~Brezzi, L.~D. Marini, and A.~Russo.
\newblock Equivalent projectors for virtual element methods.
\newblock {\em Computers $\&$ Mathematics with Applications}, 66(3):376--391,
  2013.

\bibitem{VEM-2013}
L.~Beir{\~a}o~da Veiga, F.~Brezzi, A.~Cangiani, G.~Manzini, L.~D. Marini, and
  A.~Russo.
\newblock Basic principles of virtual element methods.
\newblock {\em Mathematical Models and Methods in Applied Sciences},
  23(01):199--214, 2013.

\bibitem{VEM-2014}
L.~Beir{\~a}o~da Veiga, F.~Brezzi, L.~D. Marini, and A.~Russo.
\newblock The hitchhiker's guide to the virtual element method.
\newblock {\em Mathematical Models and Methods in Applied Sciences},
  24(08):1541--1573, 2014.

\bibitem{CVEM-elliptic}
L.~Beir{\~a}o~da Veiga, F.~Brezzi, L.~D. Marini, and A.~Russo.
\newblock Virtual element method for general second-order elliptic problems on
  polygonal meshes.
\newblock {\em Mathematical Models and Methods in Applied Sciences},
  26(04):729--750, 2016.

\bibitem{VEM-3D-MHD}
L.~Beir{\~a}o~da Veiga, F.~Dassi, G.~Manzini, and L.~Mascotto.
\newblock The virtual element method for the 3{D} resistive magnetohydrodynamic
  model.
\newblock {\em Mathematical Models and Methods in Applied Sciences},
  33(03):643--686, 2023.

\bibitem{stability-MMM}
L.~Beir{\~a}o~da Veiga, C.~Lovadina, and A.~Russo.
\newblock Stability analysis for the virtual element method.
\newblock {\em Mathematical Models and Methods in Applied Sciences},
  27(13):2557--2594, 2017.

\bibitem{DivFree-Stokes-CVEM}
L.~Beir{\~a}o~da Veiga, C.~Lovadina, and G.~Vacca.
\newblock Divergence free virtual elements for the {S}tokes problem on
  polygonal meshes.
\newblock {\em ESAIM: Mathematical Modelling and Numerical Analysis},
  51(2):509--535, 2017.

\bibitem{DivFree-NS-CVEM}
L.~Beir{\~a}o~da Veiga, C.~Lovadina, and G.~Vacca.
\newblock Virtual elements for the {N}avier--{S}tokes problem on polygonal
  meshes.
\newblock {\em SIAM Journal on Numerical Analysis}, 56(3):1210--1242, 2018.

\bibitem{FEM2008}
S.~C. Brenner.
\newblock {\em The mathematical theory of finite element methods}.
\newblock Springer, 2008.

\bibitem{CVEM-NCVEM-elliptic}
A.~Cangiani, G.~Manzini, and O.~J. Sutton.
\newblock Conforming and nonconforming virtual element methods for elliptic
  problems.
\newblock {\em IMA Journal of Numerical Analysis}, 37(3):1317--1354, 2017.

\bibitem{chen-huang-2018}
L.~Chen and J.~Huang.
\newblock Some error analysis on virtual element methods.
\newblock {\em Calcolo}, 55:1--23, 2018.

\bibitem{MHD-introduce}
P.~A. Davidson.
\newblock {\em Introduction to magnetohydrodynamics}, volume~55.
\newblock Cambridge university press, 2017.

\bibitem{Dong-CN-MHD}
X.~Dong and Y.~He.
\newblock Optimal convergence analysis of {C}rank--{N}icolson extrapolation
  scheme for the three-dimensional incompressible magnetohydrodynamics.
\newblock {\em Computers $\&$ Mathematics with Applications},
  76(11-12):2678--2700, 2018.

\bibitem{Dong-MHD-iterative}
X.~Dong, Y.~He, and Y.~Zhang.
\newblock Convergence analysis of three finite element iterative methods for
  the 2{D}/3{D} stationary incompressible magnetohydrodynamics.
\newblock {\em Computer Methods in Applied Mechanics and Engineering},
  276:287--311, 2014.

\bibitem{Dong-SAV-MHD}
X.~Dong, H.~Huang, Y.~Huang, X.~Shen, and Q.~Tang.
\newblock Optimal convergence analysis of two {RPC}-{SAV} schemes for the
  unsteady incompressible magnetohydrodynamics equations.
\newblock {\em IMA Journal of Numerical Analysis}, page drae016, 2024.

\bibitem{MHD-book}
J.-F. Gerbeau, C.~Le~Bris, and T.~Leli{\`e}vre.
\newblock {\em Mathematical methods for the magnetohydrodynamics of liquid
  metals}.
\newblock Clarendon Press, 2006.

\bibitem{MHD-existence}
M.~D. Gunzburger, A.~J. Meir, and J.~S. Peterson.
\newblock On the existence, uniqueness, and finite element approximation of
  solutions of the equations of stationary, incompressible
  magnetohydrodynamics.
\newblock {\em Mathematics of Computation}, 56(194):523--563, 1991.

\bibitem{Hou-SAV-NS}
Y.~Han, Y.~Hou, and M.~Zhang.
\newblock Analysis of divergence-free {H}$^1$ conforming {FEM} with
  {IMEX}-{SAV} scheme for the {N}avier-{S}tokes equations at high {R}eynolds
  number.
\newblock {\em Mathematics of Computation}, 92(340):557--582, 2023.

\bibitem{He-MHD-Euler}
Y.~He.
\newblock Unconditional convergence of the {E}uler semi-implicit scheme for the
  three-dimensional incompressible mhd equations.
\newblock {\em IMA Journal of Numerical Analysis}, 35(2):767--801, 2015.

\bibitem{Li-MHD-divergence-free}
R.~Hiptmair, L.~Li, S.~Mao, and W.~Zheng.
\newblock A fully divergence-free finite element method for magnetohydrodynamic
  equations.
\newblock {\em Mathematical Models and Methods in Applied Sciences},
  28(04):659--695, 2018.

\bibitem{FEM2016}
V.~John.
\newblock {\em Finite element methods for incompressible flow problems},
  volume~51.
\newblock Springer, 2016.

\bibitem{Li-SAV-NS}
X.~Li, J.~Shen, and Z.~Liu.
\newblock New {SAV}-pressure correction methods for the {N}avier-{S}tokes
  equations: stability and error analysis.
\newblock {\em Mathematics of Computation}, 91(333):141--167, 2022.

\bibitem{Li-SAV-MHD}
X.~Li, W.~Wang, and J.~Shen.
\newblock Stability and error analysis of {IMEX} {SAV} schemes for the
  magneto-hydrodynamic equations.
\newblock {\em SIAM Journal on Numerical Analysis}, 60(3):1026--1054, 2022.

\bibitem{stability-ill}
L.~Mascotto.
\newblock Ill-conditioning in the virtual element method: Stabilizations and
  bases.
\newblock {\em Numerical Methods for Partial Differential Equations},
  34(4):1258--1281, 2018.

\bibitem{interpolation-Stokes}
J.~Meng, L.~Beir{\~a}o~da Veiga, and L.~Mascotto.
\newblock Stability and interpolation properties for {S}tokes-like virtual
  element spaces.
\newblock {\em Journal of Scientific Computing}, 94(3):56, 2023.

\bibitem{VEM-2D-MHD}
S.~Naranjo-Alvarez, L.~Beir{\~a}o~da Veiga, V.~A. Bokil, F.~Dassi, V.~Gyrya,
  and G.~Manzini.
\newblock The virtual element method for a 2{D} incompressible {MHD} system.
\newblock {\em Mathematics and Computers in Simulation}, 211:301--328, 2023.

\bibitem{Shen-SAV}
J.~Shen, J.~Xu, and J.~Yang.
\newblock The scalar auxiliary variable ({SAV}) approach for gradient flows.
\newblock {\em Journal of Computational Physics}, 353:407--416, 2018.

\bibitem{Tang-MHD}
Q.~Tang and Y.~Huang.
\newblock Local and parallel finite element algorithm based on {O}seen-type
  iteration for the stationary incompressible {MHD} flow.
\newblock {\em Journal of Scientific Computing}, 70:149--174, 2017.

\bibitem{Xu-projection-MHD}
C.~Wang, J.~Wang, Z.~Xia, and L.~Xu.
\newblock Optimal error estimates of a {C}rank--{N}icolson finite element
  projection method for magnetohydrodynamic equations.
\newblock {\em ESAIM: Mathematical Modelling and Numerical Analysis},
  56(3):767--789, 2022.

\bibitem{Wei-Fealpy}
H.~Wei and Y.~Huang.
\newblock Fealpy: Finite {E}lement {A}nalysis {L}ibrary in {P}ython.
  https://github.com/weihuayi/ fealpy, Xiangtan University, 2017-2023.

\end{thebibliography}

\end{document}